\documentclass[12pt,reqno]{amsart}

\usepackage{
	amscd,
	amsfonts,
	amsmath,
	amssymb,
	cancel,
	color,
	diagbox,
	float,
	graphicx,
	hyperref,
	multicol,
	multirow,
	soul,
	ulem,
	ytableau,
}

\newtheorem{theorem}{Theorem}
\newtheorem{proposition}[theorem]{Proposition}
\newtheorem{corollary}[theorem]{Corollary}
\newtheorem{lemma}[theorem]{Lemma}

\theoremstyle{definition}
\newtheorem{definition}[theorem]{Definition}
\newtheorem{example}[theorem]{Example}
\newtheorem{remark}[theorem]{Remark}

\newcommand{\FF}{\mathbb{F}}

\def\s{\mathsf{s}}

\def\Bf{\mathbf{b}}

\def\E{\mathbf{e}}

\def\F{\mathbf{f}}
\def\Q{\mathbf{q}}
\def\S{\mathbf{s}}
\def\t{\mathbf{t}}

\def\T{\mathrm{T}}

\def\W{\mathsf{W}}
\def\bBr{\mathsf{bBr}}
\def\bBr{\mathsf{bBr}}

\def\tJ{\mathsf{tJ}}
\def\bJ{\mathsf{bJ}}

\def\Sym{\mathsf{S}}
\def\B{\mathsf{B}}
\def\C{\mathfrak{C}}
\def\D{\mathsf{D}}

\def\M{\mathsf{M}}
\def\P{\mathsf{P}}
\def\Q{\mathsf{Q}}
\def\R{\mathsf{R}}
\def\T{\mathsf{T}}
\def\J{\mathsf{J}}
\def\tJ{\mathsf{tJ}}

\def\LP{\mathsf{LP}}

\def\Br{\mathsf{Br}}

\begin{document}

\title{Brauer and Jones tied monoids}

\author[F. Aicardi]{Francesca Aicardi}
\address{Sistiana Mare 56, 34011 Trieste, Italy.}
\email{francescaicardi22@gmail.com}

\author[D. Arcis]{Diego Arcis}
\address{Facultad de Ciencias de la Salud, Universidad Aut\'onoma de Chile - Sede Talca, 5 Poniente 1670, 3460000 Talca, Chile.}
\email{diego.arcis@uautonoma.cl}

\author[J. Juyumaya]{Jes\'us Juyumaya}
\address{IMUV, Universidad de Valpara\'{\i}so\\Gran Breta\~na 1111, 2340000 Valpara\'{\i}so, Chile.}
\email{juyumaya@gmail.com}

\date{}
\thanks{The second author is part of the research group GEMA Res.180/2019 VRIP--UA and was supported, in part, by the grant Fondo Apoyo a la Investigaci\'on DIUA179-2020. The third author was supported partially by the grant FONDECYT Regular Nro.1210011. We also thanks the computational support of the National Laboratory for High Performance Computing Chile.}

\maketitle

\begin{abstract}
We introduce a ramified monoid,  attached to each Brauer--type monoid, that is, to the symmetric group, to the Jones and Brauer monoids among others. Ramified monoids correspond to a class of tied monoids which arise from knot theory and are interesting in itself. The ramified monoid attached to the symmetric group is the  Coxeter-like version of the so--called tied braid monoid. We give a presentation of the ramified monoid attached to the Brauer monoid.  Also, we introduce and studied two tied-like monoids that cannot be described as ramified monoids. However, these  monoids can also be regarded as tied versions of the Jones and Brauer monoids.
\end{abstract}

\section{Introduction}\label{sectionone}

Three pillars of quantum invariants are the Jones polynomial \cite{JoBulAMS1985}, the Homflypt polynomial \cite{FrYeHoLiMiOcBulAMS1985} and the Kauffman polynomial \cite{KaTAMS1990}, which can be  constructed, respectively, from the Temperley--Lieb algebra \cite{TeLiPRSL1971, JoBulAMS1985}, the Hecke algebra and the BMW algebra \cite{BiWeTAMS1989}. These algebras can be obtained as deformations of monoids, namely: the Jones monoid \cite{JoBulAMS1985, LaFiCA2006}, the symmetric monoid (group) and the Brauer monoid \cite{KeJSM1989,  KuMaCEJM2006} (cf. \cite{Brauer}), respectively.

Tied links are a generalization of classical links and were introduced in \cite{AiJuJKTR2016}. In the same article, a Homflypt type invariant is defined, through the  so--called bt--algebra \cite{AiJuJKTR2016, JaPoJLMS2017, MaIMRN2018, RyHaJAC2011, AiJuAMJ2020}, which turns out to be more powerful than the Homflypt polynomial on links. The bt--algebra can be considered as a tied Hecke algebra because its defining presentation, consisting of braids and tie generators, becomes into the classical Hecke algebra when taking the tie generators as the identity. Also, the bt--algebra  is understood as a deformation of the {\it tied symmetric monoid}, see Proposition \ref{PresTSn}. This monoid can also be obtained by imposing on the tied braid monoid \cite[Definition 3.1]{AiJuJKTR2016} the condition that each elementary braid is an involution.

In \cite{AiJuMathZ2018} is defined a Kauffman type invariant for tied links that is more powerful than the Kauffman polynomial  when restricting to classical links, and  also a tied BMW algebra  denoted by tBMW. This algebra is defined by a  presentation that  consists of  four types of generators: braid generators, tangle generators, tie generators and tied tangle generators. Making in this presentation the tied tangle generators as tangle generators and the tie generators as   the identity, we obtain the BMW algebra. Taking into account that the BMW can be considered as a deformation of the Brauer monoid, the task arises to define a monoid that is for the tBMW algebra as the Brauer monoid is for the BMW algebra. This task is the primary driver of the paper.

The tied braid monoid has motivated the construction of others tied braid--like monoids. Indeed, in \cite{ArJuarXiv2020},  several families  of tied monoids  were built,  which are semi direct products  of a braid--like monoid with a a suitable monoid of set partitions. However, these semi direct products do not produce a tied version of the Brauer monoid or the Jones monoid.

This article deals with the construction of tied monoids from structures that lie in the partition monoid $\C_n$  \cite{KuMaCEJM2006}. Note that $\C_n$ is the monoidal version of the  so--called partition algebra introduced indepently by V. Jones \cite{JoPro1994} and P. Martin \cite{MaJKTR1994}. More precisely, we introduce here a certain tied monoid $\R\M$, cf. \cite[p. 476]{MaEaMathZ2004}, called ramified monoid,  which is constructed from a given submonoid $\M$ of $\C_n$. We are interested in  Brauer-type monoids which have a presentation by generators and relations. Observe that, a priori, it is not a trivial matter to construct a presentation of $\R\M$ even if we know a presentation of $\M$.

When $\M$ is the symmetric group, we prove that $\R\M$ is the tied symmetric group, see Theorem \ref{isoTSRS}. For $\M$  the Brauer monoid, we obtain a presentation of $\R\M$ (Theorem \ref{PresenPBrn}), which consists of adding to the Brauer generators, the tie generators and the tied tangles generators together with certain relations from the defining presentation of the algebra tBMW \cite[Section 5]{AiJuMathZ2018}.

In the case that $\M$ is the Jones monoid, its ramified monoid does not admit a presentation by adding to the standard generators the tie generators. The reason is that the tie generators cannot be extended to non consecutive strands, since the symmetric group is absent in $\J_n$. However, one may ask what happens if we add to the Jones monoid the tie generators. This is studied in \cite{AiJuPa2019}. Now, one may add the ties and substitute the tangle generators with the tied tangle generators. The obtained monoid is denoted by $\tJ_n$ and is studied in Section \ref{sectionfiv}. It turns out that the partitions induced by the ties are easy to be described and counted. Similarly, one may consider doing the same with the Brauer monoid. i.e., by adding tie generators and substituting tangles by tied tangles. The resulting monoid is named $\bBr_n$ and it is studied in Section \ref{sectionsev}.

The paper is organized as follows. Section \ref{sectiontwo} recalls the definitions of the Jones and  Brauer monoids, as well as some properties of them. Section \ref{sectionthr} is devoted to monoids of set partitions and their relatives: the monoids of linear partitions and double set partitions. Also, we recall the partition monoid. All these monoids are essential for the development of the paper. In Section \ref{sectionfou} is introduced the main object of the paper: the ramified monoid (Definition \ref{RM}). In Section \ref{sectionfiv} it is proved that the tied symmetric monoid is the ramified monoid of the symmetric group, see  Theorem \ref{isoTSRS}. Notably, also is constructed a presentation of the ramified monoid associated to the Brauer monoid, see Theorem \ref{PresenPBrn}. In Sections \ref{sectionsix} and  \ref{sectionsev} are studied, respectively, the monoids $\bBr_n$ and $\tJ_n$.

\section{The Jones and Brauer monoids}\label{sectiontwo}

We begin the section by first giving some notations. Later, we recall the definitions of the Jones and Brauer monoids.

\subsection{}

For a monoid $\M$, we denote its identity by $1_{\M}$, but if there is no possibility of confusion, we simply write it by $1$. The group of units of $\M$ is denoted by $\M^{\times}$.

\subsection{}

Let $\J_n$ be the Jones monoid \cite{JoIM1983,LaFiCA2006}, presented by generators $\t_1,\ldots,\t_{n-1}$, called tangles, subject to the relations:
\begin{align}
\t_i^2 &=\t_i\quad\text{for all }i,\label{T1}\\
\t_i\t_j&=\t_j\t_i\quad\text{if }|i-j|>1,\label{T2}\\
\t_i\t_j\t_i&=\t_i\quad\text{if }|i-j|=1.\label{T3}
\end{align}
The size of $\J_n$ is the $n$th Catalan number $C_n=\frac{1}{n+1}\binom{2n}{n}$, see \cite[Aside 4.1.4]{JoIM1983}.

Let $\Sym_n$ be the symmetric group on $n$ symbols, which is presented as a Coxeter group, by generators $\S_1,\ldots,\S_{n-1}$ subject to the relations:
\begin{align}
\S_i^2&=1\quad\text{for all }i,\label{S1}\\
\S_i\S_j&=\S_j\S_i\quad\text{if }\vert i-j\vert>1,\label{S2}\\
\S_i\S_j\S_i&=\S_j\S_i\S_j\quad\text{if }\vert i-j\vert=1.\label{S3}
\end{align}

Let $\Br_n$ be the Brauer monoid \cite{Brauer,KuMaCEJM2006} which is defined  by  generators $\S_1,\ldots,\S_{n-1}$, $\t_1,\ldots,\t_{n-1}$ subject to the relations (\ref{S1})--(\ref{S3}) and (\ref{T1})--(\ref{T3}) together with the following relations:
\begin{align}
\t_i\S_i&=\S_i\t_i=\t_i\quad \text{for all }i,\label{Br1}\\
\t_i\S_j&=\S_j\t_i\quad\text{if }|i-j|>1,\label{Br2}\\
\S_i\t_j\t_i&=\S_j\t_i,\, \t_i\t_j\S_i=\t_i\S_j\quad\text{if }|i-j|=1.\label{Br3}
\end{align}
For each $i,j$ s.t. $|i-j|=1$, the above relations imply:
\begin{align}
\S_i\t_j\S_i&=\S_j\t_i\S_j,\label{SitjSi}\\
\t_i\S_j\t_i&=\t_i,\label{tiSjti}\\
\S_i\S_j\t_i&=\t_j\S_i\S_j=\t_j\t_i.\label{SiSjti}
\end{align}
Recall that $|\Br_n|=(2n-1)!!$ and that $\Br_n^{\times}=\Sym_n$ (see \cite[Lemma 2.6]{MaCA2002}).

\begin{figure}[H]
\includegraphics{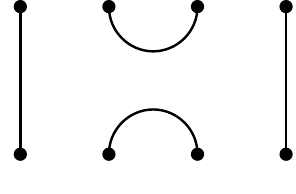}
\qquad\qquad\qquad
\includegraphics{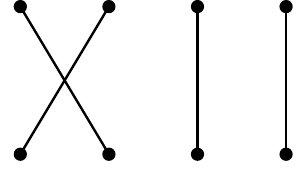}
\caption{For $n=4$, the first diagram represent $\t_2$ and the other one $\S_1$.}
\label{Fig1}
\end{figure}

The elements of $\Br_n$ can be represented by the so-called Brauer diagrams, as show the Figure \ref{Fig1}.

\section{Set partitions}\label{sectionthr}

Here we give the necessary background on set partitions, the monoid of set partitions and the partition monoid.

\subsection{}

A set partition or simply a partition of a set $X$, is a collection of nonempty subsets disjoint from each other, whose union is $X$; the members of the set partitions are called blocks. We denote by $\P(X)$ the collection of set partitions of $X$. On $\P(X)$ we have a partial order $\preceq$ defined by: $J\preceq I$ if each block of $I$ can be obtained as an union of blocks of $J$. Further, the set $\P(X)$ equipped with the product by refinement becomes an idempotent commutative monoid. This monoid will be denoted by $\P_X$ and its unity is the set partition whose blocks are singletons.

\subsection{}

For a positive integer $n$, put $[n]=\{1,\ldots,n\}$ and denote by $\P_n$ the monoid $\P_{[n]}$. The size of $\P_n$ is the $n$th Bell number $b_n$, see OEIS A000110. The elements of $\P_n$ are commonly represented as linear graphs: the vertices of the linear graph are the elements of $[n]$ and its set of edges, called arcs, are what connect $i$ with $j$, if $j$ is the minimum number in the same block of $i$ satisfying $i<j$.

\begin{figure}[H]
\includegraphics{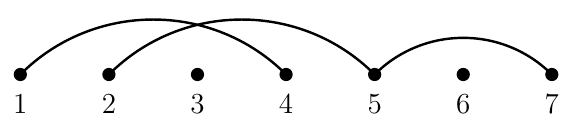}
\caption{Linear graph of $\{\{1,4\},\{2,5,7\},\{3\},\{6\}\}$.}
\end{figure}

A reference for standard facts on $\P_n$ is \cite{MaBook}.

For our purposes it is more convenient to represent the elements of $\P_n$ by a {\it diagram of ties}. A diagram of ties can be obtained by replacing the $n$ vertices and the arcs of the linear graph of a set partition by, respectively, $n$ parallel lines and dashed lines connecting the vertical lines coming from connected vertices.

\begin{figure}[H]
\includegraphics{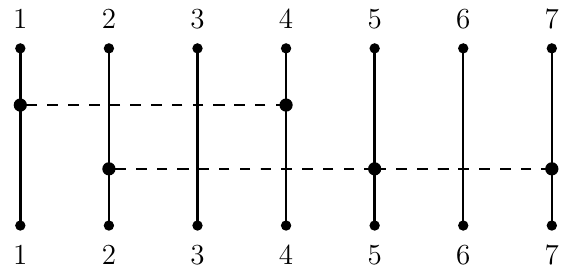}
\caption{Ties diagram of $\{\{1,4\},\{2,5,7\},\{3\},\{6\}\}$.}
\label{003}
\end{figure}

For $A\subseteq[n]$,  put $\E_A:=\{A\}\cup\{\{k\}\mid k\in[n]\backslash A\}$.

We will use the following result due to D. FitzGerald.
\begin{theorem}[{\cite[Theorem 2]{FiBAUMS2003}}]\label{Gerald}The monoid $\P_n$ can be presented with generators $\E_{i,j}:=\E_{\{i,j\}}$, where $i,j\in [n]$ with $i<j$, and the following relations:\begin{align}
\E_{i,j}^2&=\E_{i,j},\label{P1}\\
\E_{i,j}\E_{r,s}&=\E_{r,s}\E_{i,j},\label{P2}\\
\E_{i,j}\E_{j,k}&=\E_{i,j}\E_{i,k}=\E_{j,k}\E_{i,k}\qquad\text{with }i<j<k.\label{P3}
\end{align}
\end{theorem}
In what follows we will simply denote $\E_i$ instead of $\E_{i,i+1}$.

\begin{proposition}[Normal form, {\cite[Proposition 2.8]{ArJuarXiv2020}}]\label{FNPn}
For every block $K$ of an element $I$ of $\P_n$, the partition $\E_K$ has a unique decomposition\[\E_K=\E_{i_1,i_2}\cdots\E_{i_{ t-1},i_t},\]where $K=\{i_1<\cdots<i_t\}$. Thus, every set partition has a uniquely determined decomposition keeping the order of the blocks.
\end{proposition}

\subsection{}\label{subsecLP}A {\it linear partition} is a set partition of $[n]$ in which the elements of each block are consecutive. For example $\{\{1\},\{2\},\{3,4,5,6\},\{7,8\},\{9\}\}$ is a linear partition of $[9]$.

We denote by $\LP_n$ the submonoid of $\P_n$ formed by linear partitions of $[n]$.

\begin{lemma}\label{dimlin}
The number of linear partitions of $[n]$ into $k$ blocks is $\binom{n-1}{k-1}$.  Thus, the size of $\LP_n$ is $2^{n-1}$.
\end{lemma}
\begin{proof}
The proof follows by taking account that there is a bijection between linear partitions and compositions of numbers. See OEIS A000079.
\end{proof}

\begin{lemma}\label{linear}  The   monoid  $\LP_n$   is generated by  $\E_1,\ldots,\E_{n-1}$.
\end{lemma}
\begin{proof}
For a  subset $K=\{a_1,\ldots,a_{t+1}\}$ of  consecutive elements of $[n]$, Proposition \ref{FNPn} implies that $\E_K=\E_{a_1}\cdots\E_{a_t}$. Since $I=\E_{I_1}\cdots\E_{I_k}$ for all $I=\{I_1,\ldots,I_k\}\in\P_n$, then $\LP_n$ is generated by $\E_1,\ldots,\E_{n-1}$.
\end{proof}

\subsection{}

A {\it double partition}, or a {\it$2$--ramified partition}, of a set $X$ is a pair $(I,R)$ such that $I,R\in\P(X)$ and $I\preceq R$, cf. \cite[p. 476]{MaEaMathZ2004}. We denote by $\D\P(X)$ the collection of double partitions of $X$.

\begin{remark}\label{RamDoub}
Observe that if $(I,R)$ is a double partition of $X$, then $R$ determines a unique set partition $R(I)$ of $I$, in which two blocks of $I$ belong to the same block of $R(I)$ if they are  contained in the same block of $R$. Moreover, there is a bijection between double partitions of $X$ and the set of pairs $(I,R(I))$ with $I\in\P(X)$ and $R(I)\in\P(I)$.
\end{remark}

The set $\D\P(X)$ equipped with the product inherited of $\P_X\times\P_X$ becomes an idempotent commutative monoid with unity $(1,1)$. This monoid will be denoted by $\D\P_X$.

\subsection{}

For a positive integer $n$, denote by $\D\P_n$ the monoid $\D\P_{[n]}$.

Note that for each $I\in\P_n$ with $k=|I|$, there are $b_k$ set partitions $R(I)$ of $I$ such that $(I,R)\in\D\P_n$. It is well known that the number of set partitions of $[n]$ into $k$ blocks corresponds to the Stirling number $S(n,k)$. Hence, the size of $\D\P_n$ is $\sum_{k=1}^nS(n,k)b_k$, see OEIS A000258.

\begin{remark}\label{DoubNorm}
Since $(I,J)=(I,I)(1,J)$ for all $(I,J)\in\D\P_n$, then $\D\P_n$ is generated by the pairs $(\E_{i,j},\E_{i,j})$ and $(1,\E_{r,s})$, where $i,j,r,s\in[n]$ with $i<j$ and $r<s$. Thus, the normal form of $\P_n$ (Proposition \ref{FNPn}) induces a natural normal form on $\D\P_n$. Indeed, each $(I,J)\in\D\P_n$ is represented by a word $uv$, where $u$ inherits the normal form of $I$ in the generators $(\E_{i,j},\E_{i,j})$ and $v$ inherits the normal form of $J$ in the generators $(1,\E_{i,j})$.
\end{remark}

\begin{proposition}
The monoid $\D\P_n$ has a presentation with generators $a_{i,j}:=(e_{i,j},e_{i,j})$ and $b_{r,s}:=(1,e_{r,s})$, where $i,j,r,s\in[n]$ with $i<j$ and $r<s$, both satisfying relations {\normalfont(\ref{P1})} to {\normalfont(\ref{P3})}, and subject to the following additional relations:\begin{equation}\label{EqPDoubP}a_{i,j}b_{r,s}=b_{r,s}a_{i,j},\qquad a_{i,j}b_{i,j}=a_{i,j}.\end{equation}
\end{proposition}
\begin{proof}
Remark \ref{DoubNorm} implies that each double partition of $[n]$ can be written as normal form by using only relations {\normalfont(\ref{P1})} to {\normalfont(\ref{P3})} together with {\normalfont(\ref{EqPDoubP})}.
\end{proof}

\subsection{\it The partition monoid}

Set $[n']=\{1',\ldots,n'\}$ where $k':=n+k$ for all $k\in[n]$. The {\it partition monoid} $\C_n$ is the set $\P([n]\cup[n'])$ equipped with the usual product by {\it concatenation} $\ast$. A formal definion of this product is given in \cite[Definition 2.2]{XiC1999M}. Note that the identity of $\C_n$ is the partition formed by the blocks $\{i,i'\}$, where $1\leq i\leq n$. For more information on this monoid, see for example \cite{MaCA2002,KuMaCEJM2006}.

If a block of a  set partition of $\C_n$ has cardinal $2$, it is called an {\it arc}. An arc is called a {\it bracket} if it is contained either in $[n]$ or in $[n']$, otherwise it is called a line, i.e. an arc that intersects both $[n]$ and $[n']$.


\begin{remark}\label{SJBsubmonoidCn}
It is well known that the symmetric group and the Jones and Brauer monoids can be conceived as submonoids of the partition monoid, see for instance \cite{KuMaCEJM2006}. Indeed:
\begin{enumerate}
\item$\Br_n$ can be regarded as the submonoid $\mathfrak{B}_n$ of $\C_n$ formed by the partitions whose blocks are only arcs. See Lemma \ref{BLBr} for details.
\item$\J_n$ is the submonoid of $\mathfrak{B}_n$ whose elements are known as planar partitions.
\item$\Sym_n$ is the submonoid of $\mathfrak{B}_n$ formed by partitions whose blocks are  only lines. It is well known that $\C_n^{\times}=\Sym_n$, see \cite{MaMaMB2007}.
\end{enumerate}
\end{remark}

\begin{remark}\label{CompPre}
Note that the elements of $\Br_n$ are not comparable through $\preceq$, indeed if $I\in\Br_n$ and $I\prec J$, for some $J\in P([n]\cup[n'])$, there is $B\in J$ with $|B|\geq4$. Thus $J\not\in\Br_n$.
\end{remark}

\subsection{\it Diagrammatic representation }

\begin{figure}[H]
\includegraphics{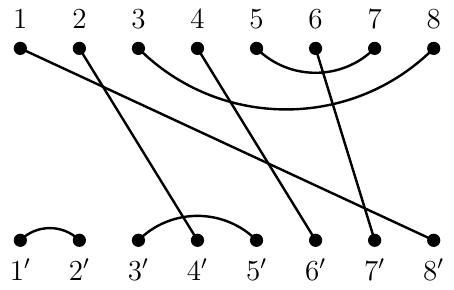}
\qquad
\includegraphics{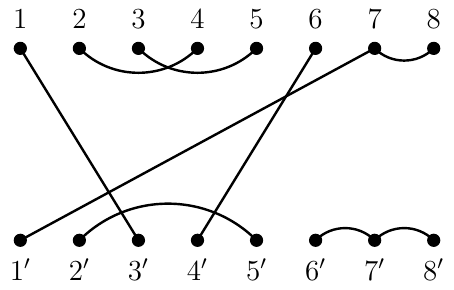}
\caption{ Two diagrams representing elements of the partition monoid.}
\label{PartitionsCn}
\end{figure}

It is convenient to represent the set partitions of $\C_n$ by placing the elements of $[n]\cup[n']$ in two parallel lines, one containing the elements of $[n]$ and the other one containing those of $[n']$. Thus, blocks of set partitions in $\C_n$ that intersect both $[n]$ and $[n']$ can be represented as two usual blocks of $[n]$ and $[n']$, joined by at least one line connecting an element of $[n]$ with one of the $[n']$. Here we will choose this line as the leftmost one, see Figure \ref{PartitionsCn}, however it can also be chosen, for example, as the rightmost one or even both lines may be considered, see Figure \ref{tiesmapFig}.

With this diagrammatic representation, the product by concatenation $D\ast D'$ can be obtained by the following procedure: we place $D$ above $D'$ and so identifying $i'$ of $D$ with $i$ of $D'$, for all $i\in [n]$.
The diagram of $D\ast D'$ has the $n$ elements $[n]$ of the top line of $D$ and the $n$ elements $[n']$ of the bottom line of $D'$. Two of such elements in $[n] \cup [n']$ belong to the same block of $D\ast D'$ if they are connected by a sequence of lines from the diagrams of $D$ or of $D'$. The representation of the product $D\ast D'$ is the diagram obtained by removing all loops that appear in the middle row.

\begin{figure}[H]
\[\vcenter{\hbox{
\includegraphics{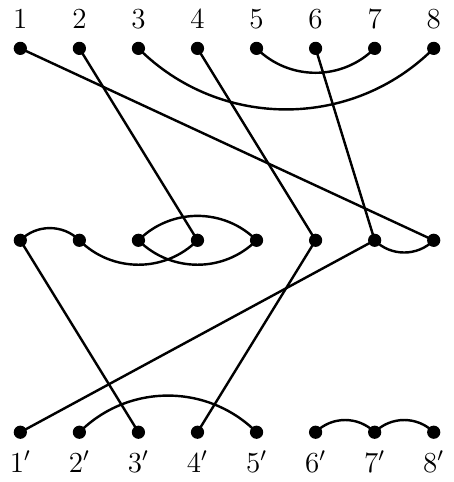}
}}
\qquad=\quad
\vcenter{\hbox{
\includegraphics{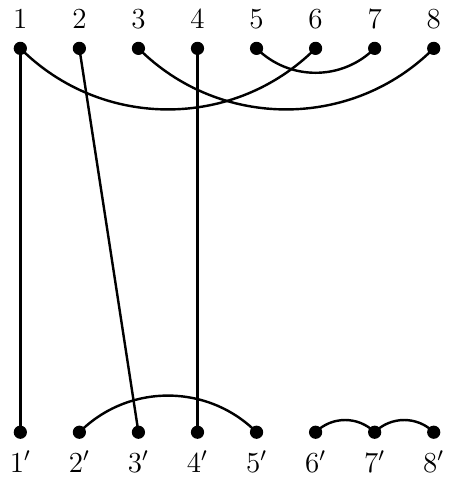}
}}
\]\caption{Concatenation product of the set partitions represented by the diagrams in Figure \ref{PartitionsCn}.}
\end{figure}

\section{Ramified monoids}\label{sectionfou}

In this section we introduce the ramified monoid $\R\M$ associated to a submonoid $\M$ of $\C_n$. Diagramatically, the product of this new monoid can be understood by expressing the  concatenation  product $\ast$ in terms of the product by refinement, see (\ref{circledastast}). It is worth noting that (theoretically) the ramified monoid construction can be applied to any finite group or finite monoid. This is a consequence of the Cayley theorem and its generalization for monoids by observing that both the symmetric group and the full transformation monoid can be considered as submonoid of the partition monoid.

\subsection{}

Let $A=\{a_1,\ldots,a_n\}$ be a set which is disjoint with $[2n]$. We set two maps $f^A, f_A$ from $[2n]$ to $[2n]\cup A$ as follows.
\[f^A(k)=\left\{ \begin{array}{ll}
a_k&\text{if $1\leq k\leq n$}\\
k&\text{otherwise}\end{array}\right.\text{and}\quad
f_A(k)=\left\{\begin{array}{ll}
a_{n-k}&\text{if $n+1\leq k\leq 2n$}\\
k&\text{otherwise.}\end{array}\right.\]
For $I\in\P([2n])$, we set:
\[I^A=\{f^A(B)\vert B\in I\}\quad\text{and}\quad I_A=\{f_A(B)\vert B\in I\}.\]
Note that $I^A\in P(A\cup[n+1,2n])$ and $I_A\in P([n]\cup A)$.   With this data we obtain the following relation\begin{equation}\label{circledastast}
I\ast J=(I_A\cdot J^A)\cap[2n],\quad\text{for all }I,J\in\C_n.
\end{equation} where the product $I_A\cdot J^A$ occurs in $P_{A\cup[2n]}$ by considering $I_A$ and $J^A$ as set partitions of $A\cup[2n]$.

\begin{definition}[Cf. \cite{MaEaMathZ2004}]\label{RM}
Let $\M$ be a submonoid of $\C_n$. The {\it ramified monoid} of $\M$, denoted by $\R\M$, is the set formed by the pairs $(I,R)\in\D\P([2n])$ with $I\in\M$, endowed with the product below.\[(I,R)(J, S)=(I\ast J,R\ast S),\]where $(J,S)\in\D\P([2n])$ and $J\in\M$.
\end{definition}

Notice that $\M$ embeds in $\R\M$ via the map $I\mapsto(I,I)$ for all $I\in\M$. In particular, we identify $\Sym_n$ with the set of pairs $(s,s)\in\R\C_n$ such that $s\in\Sym_n$.

\begin{proposition}\label{PropiPM}
For every submonoid $\M$ of $\C_n$, we have:
\begin{enumerate}
\item $(\R\M)^{\times}\cong\M^{\times}$.
\item The size of $\R\M$ is $|\M|b_n$.
\end{enumerate}
\end{proposition}
\begin{proof}
Suppose that $(I,R)\in(\R\M)^{\times}$, then there exist $(J, S)$ such that $(I,R)(J, S)=(1,1)$. This implies that $I\in\M^{\times}\subseteq\Sym_n$ and $R\in\C_n^{\times}=\Sym_n$.  Since $I\preceq R$ and $I,R\in \Sym_n$, then $I=R$. So, $(\R\M)^{\times}$ is formed by the pairs $(I,I)$ with $I\in\M^{\times}$. Thus, $(\R\M)^{\times}$ and $\M^{\times}$ are naturally isomorphic.

The claim (2) follows from Remark \ref{RamDoub}.
\end{proof}
\begin{definition}\label{DefEij}
For $0<i,j<n$, let $E_{i,j}$ be the partition of $[n]\cup[n']$ whose blocks are $\{i,j,i',j'\}$ and $\{k,k'\}$ for $k\not=i$. Further, define
\begin{equation}\label{TildeLiE_i}
\widetilde{E}_{i,j}=(1,E_{i,j}).
\end{equation}
\end{definition}
Observe that $E_{i,i}=1$ and $E_{i,j}=E_{j,i}$. In the case $\vert i-j\vert=1$, we shall denote $E_k$ instead of $E_{i,j}$, where  $k:=\mathrm{min}\{i,j\}$.

\begin{proposition}\label{SnPnPM}
For every submonoid $\M$ of $\C_n$, the mapping $\E_{i,j}\mapsto \widetilde{E}_{i,j}$ defines a monomorphism of $\P_n$ in $\R\M$.
\end{proposition}
\begin{proof}
The proof follows by observing that the natural bijection between $[n]$ and $\{\{i,i'\}|1\leq i\leq n\}$ is compatible with the union of sets.
\end{proof}

\begin{figure}[H]
\[\vcenter{\hbox{
\includegraphics{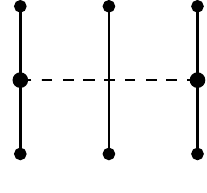}
}}
\quad\to\quad
\left(\vcenter{\hbox{
\includegraphics{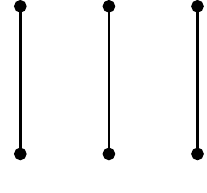}\,\,\,
}}
,
\vcenter{\hbox{
\includegraphics{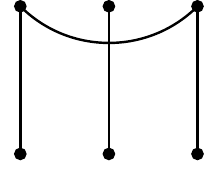}\,\,
}}\right)\]
\caption{The image of $\E_{1,3}$ under the monomorphism of Proposition \ref{SnPnPM}.}
\label{tiesmapFig}
\end{figure}

\section{Presentations of two ramified monoids}\label{sectionfiv}

In this section we consider the ramified monoid $\R\M$ for $\M$ the symmetric group and the Brauer monoid. More precisely, firstly we prove that $\R\Sym_n$ coincides with the tied monoid of the symmetric group as  constructed in \cite{ArJuarXiv2020}. Secondly, we give a presentation of the monoid $\R\Br_n$. This presentation is obtained by showing that $\R\Br_n$ is isomorphic to an abstract monoid $\Q_n$ whose relations in its defining presentation come from the tBMW algebra.

\subsection{\it The ramified monoid  $\R\Sym_n$}

\begin{definition}
For $0<i<n$, let $L_i$ be the partition consisting of blocks: $\{i,i'+1\}$, $\{i+1,i'\}$ and $\{j,j'\}$ for $j \not= i,i+1$. Moreover, define the double partitions\[\widetilde{L}_i=(L_i,L_i).\]
\end{definition}

Note that $\Sym_n$ is generated by the $L_i$'s. Indeed $L_i$ corresponds to the generator $\S_i$ of $\Sym_n$, see (3)Remark \ref{SJBsubmonoidCn}.
Further, observe that\begin{equation}\label{LiEj}
L_i\ast E_{j,k}=E_{\S_i(j),\S_i(k)}\ast L_i,
\end{equation}where $\S_i(j)$ is the effect of the transposition $(i,i+1)$ on $j$.

\begin{proposition}\label{GeneraPSn}
$\R\Sym_n$ is generated by the $\widetilde{L}_i$'s and the $\widetilde{E}_i$'s, where  $\widetilde{E}_i=(1,E_i)$.
\end{proposition}
\begin{proof}
By using (\ref{LiEj}) we have: \begin{equation}\label{LS}\widetilde{L}_i\widetilde{E}_{j,k}=(L_i\ast 1, L_i\ast E_{j,k})=(1\ast L_i,E_{\S_i(j),\S_i(k)}\ast L_i).\end{equation}
Therefore, for every $i,j$ such that $j>i+1$, we get
\begin{eqnarray*}
\widetilde{E}_{i,j}&=&(L_{j-1},L_{j-1})
\cdots(L_{i+1},L_{i+1})(1,E_i)
(L_{i+1},L_{i+1})
\cdots(L_{j-1},L_{j-1})\\
&=&(1,L_{j-1}\ast\cdots\ast L_{i+1}\ast E_i\ast L_{i+1}\ast\cdots\ast L_{j-1}).
\end{eqnarray*}
The proof thus follows from Proposition \ref{SnPnPM}.
\end{proof}

\subsubsection{}

The {\it tied symmetric monoid} $\T\Sym_n$ is the semi direct product between $\Sym_n$ and $\P_n$, that defines the permutation action of $\Sym_n$ on $\P_n$, that is, the action induced from the natural action of $\Sym_n$ on $[n]$. The Lavers' method \cite{La1998} applied to the Coxeter presentation of $\Sym_n$ and the FitzGerald presentation of $\P_n$ yields the following proposition.

\begin{proposition}[Cf. {\cite[Proposition 4.8]{ArJuarXiv2020}}]\label{PresTSn}
$\T\Sym_n$ is presented by generators $\S_1,\S_2,\ldots,\S_{n-1}$, $\E_1, \ldots,\E_{n-1}$ subject to the relations (\ref{S1})--(\ref{S3}), (\ref{P1}), (\ref{P2}) and:  
\begin{align}
\E_i^2& =\E_i \quad \text{for all }i,\label{Ei2}\\
\E_i\E_j&=\E_j\E_i\quad  \text{for all $i$ and $j$}, \label{EiEj}\\
\E_i\S_j\S_i & =\S_j\S_i\E_j\quad\text{if $|i-j|=1$,}\label{TSn1}\\
\S_i\E_j & =\E_j\S_i\quad\text{if $|i-j|\neq 1$,}\label{TSn2}\\
\E_i\E_j\S_i&=\E_j\S_i\E_j=\S_i\E_i\E_j\quad\text{if $|i-j|=1$.}\label{TSn3}
\end{align}
\end{proposition}
Observe that the monoid algebra of $\T\Sym_n$ is the specialization $\mathcal{E}_n(1)$ of the bt--algebra, see \cite{BaART2013}.

\begin{remark}\label{TBn}\rm
Observe that $\T\Sym_n$ can be also defined as the quotient obtained from the tied braid monoid $\T\B_n$ \cite{AiJuJKTR2016, ArJuarXiv2020} under the congruence generated by $\sigma_i=\sigma_i^{-1}$ for every elementary braid $\sigma_i$ of the braid group on $n$ strands.
\end{remark}

\begin{theorem}[Cf. {\cite[Theorem 4.2]{BaART2013}}]\label{isoTSRS} The map $\phi:\S_i\mapsto\widetilde{L}_i $, $\E_i\mapsto \widetilde{E}_i$, defines a monoid isomorphism from $\T\Sym_n$ to $\R\Sym_n$.
\end{theorem}
\begin{proof}
By using (\ref{LiEj}), it is a routine to check that the map  preserves with the defining relations of $\T\Sym_n$. Proposition \ref{GeneraPSn} implies that $\phi$ is an epimorphism. We will see that $\phi$ is injective. Because $\T\Sym_n=\Sym_n\ltimes\P_n$, every element of  $\T\Sym_n$ can be written as $(r,e)$, with $r\in\Sym_n$ and $e\in\P_n$. Suppose $\phi(r,e)=\phi(s,f)$, where $(r,e),(s,f)\in\T\Sym_n$. We have $\phi(r,e)=\phi(r,1_{\P_n})\phi(1_{\Sym_n},e)$ and $\phi(s,f)=\phi(s,1_{\P_n})\phi(1_{\Sym_n},f)$. Therefore $r=s$. Now, as $r\in\Sym_n$ it follows that $e=f$. Thus the proof is concluded.
\end{proof}

\subsection{\it The ramified monoid $\R\Br_n$ }

To show a presentation of $\R\Br_n$ we need to introduce first an auxiliary monoid $\Q_n$.
\begin{definition}
Let $n\geq2$, we define $\Q_n$ as the monoid presented with generators $\S_1,\ldots,\S_{n-1}$, $\t_1,\ldots,\t_{n-1}$, $\E_1,\ldots,\E_{n-1}$, $\F_1,\ldots,\F_{n-1}$, satisfying the relations (\ref{T1})--(\ref{Br3}) and (\ref{Ei2})--(\ref{TSn3}), together with the following relations:
\begin{align}
\F_i^2&=\F_i&&\text{for all }i,\label{Fi2}\\
\F_i\F_j&=\F_j\F_i&&\text{for all $i$ and $j$ s.t.}  \vert i-j\vert >1,\label{FiFj}\\
\E_i\F_i&=\F_i\E_i=\F_i&&\text{for all }i,\label{EiFi}\\
\E_i\F_j&=\F_j\E_i&&\text{for all }i,j,\label{EiFj}\\
\F_i\F_j\F_i&=\E_j\F_i\E_j&&\text{for all $i$ and $j$ s.t. }\vert i-j\vert=1,\label{FiFjFi}\\
\E_i\t_i&=\t_i\E_i=\t_i&&\text{for all }i,\label{Eiti}\\
\E_i\t_j&=\t_j\E_i&&\text{for all $i$ and $j$ s.t. }\vert i-j\vert>1,\label{Eitj}\\
\F_i\t_j&=\t_j\F_i&&\text{for all $i$ and $j$ s.t. }\vert i-j\vert>1,\label{Fitj}\\
\F_i\S_j&=\S_j\F_i&&\text{for all $i$ and $j$ s.t. }\vert i-j\vert>1,\label{FiSj}\\
\F_i\E_j&= \E_j\F_i =\E_j\t_i\E_j&&\text{for all $i$ and $j$ s.t. }\vert i-j\vert=1,\label{FjEi}\\
\S_i\F_i&=\F_i\S_i=\F_i&&\text{for all }i,\label{SiFi}\\
\S_i\F_j\S_i&=\S_j\F_i\S_j&&\text{for all $i$ and $j$ s.t. }\vert i-j\vert=1,\label{SiFjSi}\\
\F_i\t_i&=\t_i\F_i=\t_i&&\text{for all }i.\label{Fiti}
\end{align}
\end{definition}

For all $i$ and $j$ s.t. $\vert i-j\vert=1$, the defining relations of $\Q_n$ imply:
\begin{align}
\F_i\F_j&=\E_j\t_i\t_j\E_i,\label{FiFj}\\
\F_i\t_j&=\E_j\t_i\t_j,\quad\t_i\F_j=\t_i\t_j\E_i,\label{Fitj}\\
\E_i\t_j &=\F_j\S_i\t_j,\quad \t_i\E_j=\t_i\S_j\F_i,\label{Eitj}\\
\F_i\S_j\F_i&=\E_j\t_i\E_j,\label{FiSjFi}\\
\S_i\E_j\S_i&=\S_j\E_i\S_j,\label{SiEjSi}\\
\F_i\F_j\E_i&=\F_i\F_j ,\label{FiFjEi}\\
\E_i\F_j\F_i&=\F_j\F_i ,\label{EiFjFi}\\
\F_i\S_j\S_i&=\S_j\S_i\F_j.\label{SiSjFi}
\end{align}

\begin{proposition}\label{submon}
$\T\Sym_n$ and $\Br_n$ are submonoids of $\Q_n$. In consequence, $\P_n$ is also a submonoid of $\Q_n$.
\end{proposition}
\begin{proof}
Notice that  the generators $\S_i$'s and $\E_j$'s together with relations (\ref{S1})--(\ref{S3}), (\ref{Ei2}), (\ref{EiEj}) and (\ref{TSn1})-(\ref{TSn3}) of $\Q_n$, form the defining presentation of $\T\Sym_n$. Since at least one of the generators $\t_i$'s, $\F_j$'s occur on both sides of all other defining relations of $\Q_n$, it follows that the submonoid of $\Q_n$ generated by $\S_i$'s and $\E_j$'s coincides with $\T\Sym_n$. The proof that $\Br_n$ is a submonoid of $\Q_n$ is similar.
\end{proof}
 
\begin{definition}\label{extended}
For $0<i<j<n$, we call {\it extended ties} or simply {\it ties}, the following elements of $\Q_n$ defined recursively by:\begin{equation}\label{EIJ}\E_{i,i+1}:=\E_i,\quad\quad\E_{i,j}:=\s_{j-1}\E_{i,j-1}\s_{j-1},\quad j>i+1.\end{equation}
\end{definition}

The following lemma collects some properties of the extended ties.
\begin{lemma}\label{ties}
For $0<i<j<n$, we have:
\begin{enumerate}
\item[(i)]The $\E_{i,j}$ satisfy the relations (\ref{P1})--(\ref{P3}).
\item[(ii)]$\E_{i,j}=\S_{i}\E_{i+1,j}\S_i,\quad j>i+1$.
\item[(iii)]$\S_k\E_{i,j}=\E_{\s_k(i),\s_k(j)}\S_k$ and $\E_{i,j}\S_k=\S_k\E_{\s_k(i),\s_k(j)}$.
\item[(iv)]$\E_{i,j}$ commutes with $\F_k$ for every $k$.
\item[(v)]$\E_{i,j}$ commutes with $\t_k$ for every $k\notin\{i-1,i,j-1,j\}$.
\end{enumerate}
\end{lemma}
\begin{proof}
Note that the generators $\S_i$'s and $\E_i$'s of $\Q_n$ satisfy all relations of $\T\Sym_n$, and so the relations of the  well known tied braid monoid $\T\B_n$, see Remark \ref{TBn}. Hence, \cite[Lemma 2]{AiJuAiM2021} implies (i), (ii) and (iii). 

To prove (iv), it is suffices to consider the non trivial cases $i\le k\le j-1$. If $i=k$, by (\ref{EiFi}) write $\F_i=\F_i\E_i$. Hence, by using (\ref{P3}), $\F_i\E_{i,j}=\F_i\E_i\E_{i+1,j}$. Writing $\E_{i+1,j}$ in terms of $\E_{i+1}$ as in (\ref{EIJ}), the commutativity follows from (\ref{EiEj}). If $k=j-1$ write similarly $\F_{j-1}=\F_{j-1}\E_{j-1}$, $\F_{j-1}\E_{i,j}=\F_i\E_{j-1}\E_{i,j-1}$, and take the expression of $\E_{i,j-1}$ in terms of $\E_{j-2}$, by using (ii). If $i<k<j-1$, use expression (\ref{EIJ}) for $\E_{i,j}$ and then relation (\ref{SiSjFi}) twice.

For (v),the non trivial cases are when $i<k<j-1$. Proceed as at point (iv) using two times relation (\ref{SiSjti}).
\end{proof}

\subsection{}

The remainder of the section aims to give a presentation of $\R\Br_n$. This will be done in Theorem \ref{PresenPBrn} by proving that $\Q_n$ and $\R\Br_n$ are isomorphic. To state this theorem we need to introduce the following lemmas and notations.

\begin{definition}
For $0<i<n$, let $H_i$ be the partition of $\C_n$ with blocks $H_{i,k}$'s, where $H_{i,i}=\{i,i+1\}$, $H_{i,i'}=\{i',i'+1\}$ and $H_{i,k}=\{k,k'\}$ whenever $k\neq i,i'$. Moreover, define\begin{equation}\label{TildeHiF_i}\widetilde{H}_i=(H_i,H_i),\qquad\widetilde{F}_i=(H_i,E_i).
\end{equation}
\end{definition}

\begin{lemma}[{\cite[Theorem 3.1]{KuMaCEJM2006}}]\label{BLBr}
The map sending $\S_i$ to $L_i$ and $\t_i$ to $H_i$ extends to an isomorphisms  between $\Br_n$ and $\mathfrak{B}_n$.
\end{lemma}

\begin{remark}
In virtue of Lemma \ref{BLBr}, we will denote by $s\in\Sym_n$ a product of $\s_i$'s as well as the corresponding concatenation product of $L_i$'s.
\end{remark}

We will use the following normal form for $\Br_n$.
\begin{proposition}\label{QNormBrn}
For every element $g\in\Br_n$, there are uniquely defined $s,s'\in\Sym_n$ and  a unique $k\leq\frac{n}{2}$ such that
\begin{equation}\label{NFBr}
g=s\ast H_1\ast H_3\ast\cdots\ast H_{2k-1}\ast s'.
\end{equation}
\end{proposition}
\begin{proof}
Let $I$ be a set partition of $[n]$ (or of $[n']$) whose blocks have at most two elements, and let $k\leq\frac{n}{2}$ be the number of  nontrivial blocks in $I$. Note that these set partitions can be obtained by intersecting an element of $\Br_n$ with $[n]$ (or with $[n']$).
\begin{figure}[H]
\includegraphics{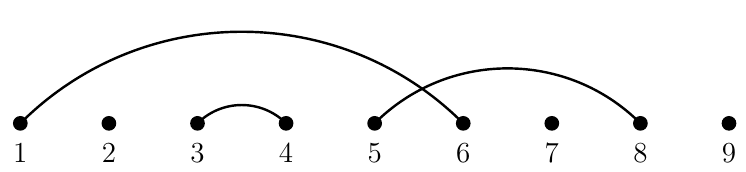}
\caption{Linear graph of $\{\{1,6\},\{2\},\{3,4\},\{5,8\},\{7\},\{9\}\}$.}
\end{figure}
Denote by $\{a_1,b_1\},\ldots,\{a_k,b_k\}$ the nontrivial blocks of $I$ satisfying $a_1<\cdots<a_k$ and $a_i<b_i$ for all $i\in[k]$, and by $\{c_1\},\ldots,\{c_l\}$ the singleton blocks of $I$, satisfying $c_j<c_{j+1}$.  We will denote by $n_I$ the unique permutation in $\Sym_n$ be the permutation satisfying:\[ n_I(a_i)=2i-1,\quad  n_I(b_i)=2i\quad\text{for}\ 1\leq i\leq k,\quad\text{and}\quad  n_I(c_j)=2k+j,\quad\text{for}\ 1\le j\leq n-2k.\] Note that $ n_I(I)$  is the unique set partition of $[n]$ with blocks of at most two elements such that $ n_I(I)\cap[2k]$ is nonnesting and noncrossing  with no trivial blocks and also $ n_I(I)\cap[2k+1,n]$ is trivial.
\begin{figure}[H]
\includegraphics{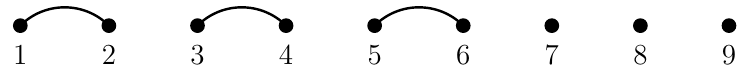}
\caption{Linear graph of $n_I(I)$ with $ I=\{\{1,6\},\{2\},\{3,4\},\{5,8\},\{7\},\{9\}\}$.}
\end{figure}

Recall that every set partition of $[n]$ can be regarded inside the partition monoid $\C_n$ by adding $n$ singleton blocks. We have $I= n_I^{-1}( n_I(I))= n_I( n_I^{-1}(I))$.\begin{figure}[H]
\includegraphics{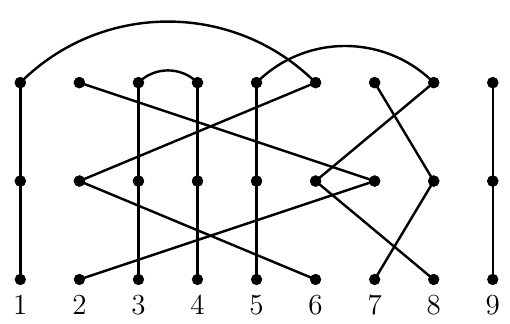}
\quad
\includegraphics{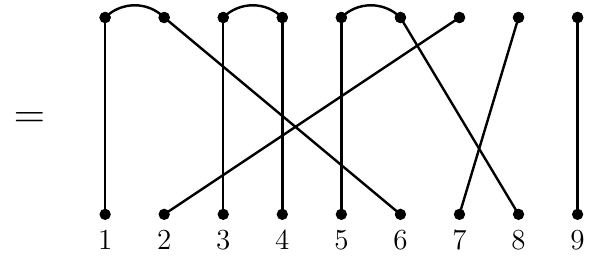}
\caption{Linear graph of $n_I^{-1}(n_I(\{\{1,6\},\{2\},\{3,4\},\{5,8\},\{7\},\{9\}\}))$.}
\end{figure}

By considering the permutations $n_I,n_J$ for $I=g\cap[n]$ and $J=g\cap[n']$, we obtain that $g=s\ast H_1\ast H_3\ast\cdots\ast H_{2k-1}\ast s'$, where $2k$ is the number of brackets of $g$, $s=\eta_I$ and $s'=t_g*\eta_J ^{-1}$ with $t_g$ the unique permutation, trivial in $[2k]$, such that the lines of $g$ are obtained by composing it with the permutations $n_I$ and $n_J^{-1}$. For instance for $g=\{\{1,5\},\{2,3\},\{4,3'\},\{6,2'\},\{1',5'\},\{4',6'\}\}\in\mathfrak{B}_6$, we have $k=2$  and the following permutations:\[n_I=
\left(\begin{smallmatrix}1&2&3&4&5&6\\1&3&4&5&2&6\end{smallmatrix}\right),\quad
t_g=\left(\begin{smallmatrix}1&2&3&4&5&6\\1&2&3&4&6&5\end{smallmatrix}\right),\quad
n_J=\left(\begin{smallmatrix}1&2&3&4&5&6\\1&5&6&3&2&4\end{smallmatrix}\right).\]
See Figure \ref{BrNorFormPic}.
\begin{figure}[H]
\[\vcenter{\hbox{
\includegraphics{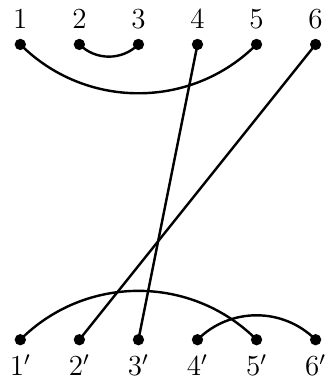}
}}\,\,=\,\,\vcenter{\hbox{
\includegraphics{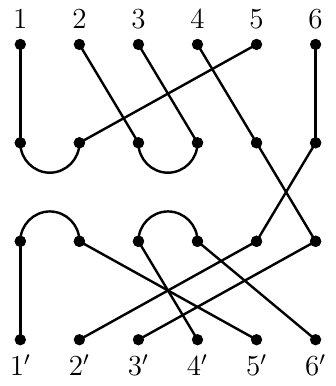}
}}\,\,=\,\,\vcenter{\hbox{
\includegraphics{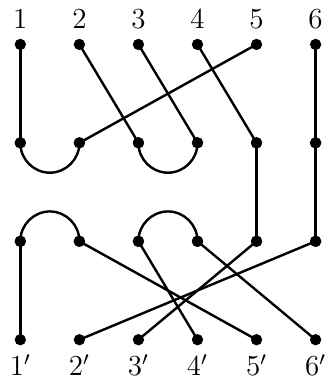}
}}\]
\caption{The set partition $\{\{1,5\},\{2,3\},\{4,3'\},\{6,2'\},\{1',5'\},\{4',6'\}\}$ can be written as $s\ast H_1\ast H_3\ast s'$  with $s,s'\in\Sym_n$.}
\label{BrNorFormPic}
\end{figure}
\end{proof}

We are going to show that $\R\Br_n$ is generated by $\widetilde{L}_i$, $\widetilde{H}_i$, $\widetilde{E}_i$ and $\widetilde{F}_i$, $i=1,\dots,n-1$. In fact, we prove the following.
\begin{lemma}\label{GenPBrn}
Every element of $\R\Br_n$ can be written as the product:\[r*T_1*T_3*\cdots*T_{2k-1}*r',\quad\text{where}\quad r,r'\in\R\Sym_n\quad\text{and}\quad T_{2i-1}\in\{\widetilde{H}_{2i-1},\widetilde{F}_{2i-1}\}.\]
\end{lemma}
\begin{proof}
Let $(I,R)\in \R\Br_n$.  Proposition \ref{QNormBrn} implies that $I=s\ast H_1\ast\cdots\ast H_{2k-1}\ast s'$.  Let $U:=\{\{a_1,b_1\},\ldots,\{a_k,b_k\}\}$ be the set of nontrivial blocks of $I\cap[n]$ satisfying $a_1<\cdots<a_k$ and $a_i<b_i$ for all $1\leq i\leq k$, and similarly $V:=\{\{a'_1,b'_1\},\ldots,\{a'_k,b'_k\}\}$ be the set of nontrivial blocks of $I\cap[n']$ satisfying $a'_1<\cdots<a'_k$ and $a'_i<b'_i$ for all $1\leq i\leq k$. Let $C:=\{\{c_j,c'_j\}\}_{j=2k+1}^n$ be the set of lines in $I$.

Let $J=I*s'^{-1}$. Starting from $J_0= J$, $G_0= J$, we will define two sequences $J_k,G_k$ ending with $J_t=I,G_t=R$. The running index $k$ increases every time we consider a new block of $I$. Denote $B_i$, $1\leq i\leq k$, the blocks of $R$ such that $B_i\cap U\not=\emptyset$, and $i$ is the minimum index of the pairs $\{a_j,b_j\}$ in them. For every block $B_i$ do the following. For every line in $\{c_j,c'_j\}\in B_i\cap C$, set $J_k=J_{k-1}$ and $G_k=E_{p,q}\ast G_{k-1}$, where $p=a_i$ and $q=c_j$. For every bracket $\{a_j,b_j\}$ in $B\cap U$, $j\not=i$, define $J_k=J_{k-1}$ and $G_{k}=E_{p,q}\ast G_{k-1}$ where $p=a_i$ or $q=a_j$. If $m$ is the minimum index such that $\{a'_m,b'_m\}\in B_i\cap V$, then $J_{k }=J_{k-1}\ast s_i^{-1}$,  where $s_i(2i-1)=a'_m$ and $s_i(2i)=b'_m$, while $G_{k }=G_{k-1} \ast E_{2i-1}\ast s_i^{-1}$. Since $H_i\ast E_i=E_i$, the generator $H_i$ in $G_{k}$ is replaced by $E_i$. For all other $\{a'_j,b'_j\}\in B_i\cap V$, $J_k=J_{k-1}$ and $G_{k}=E_{p,q}\ast G_{k-1}$, where $p=a'_m-n$, $q=a'_j-n$. Now, let $B$ any of the remaining blocks in $R$, satisfying $B\cap U=\emptyset$. If also $B\cap V=\emptyset$, then let $m$ be the minimum index such that $\{c_m,c'_m\}\in B\cap C$. Then, for every other line $\{c_j,c'_j\}\in B\cap C$, $j\not=m$, $J_k=J_{k-1}$ and $G_{k}=E_{p,q}\ast G_{k-1}$, where $p=c_m$ and $q=c_j$. Finally, if $B\cap V\not=\emptyset$, then let $m$ be the minimum index such then $\{a'_m,b'_m\}\in B\cap V$. Then, for every line in $\{c_j,c'_j\}\in B\cap C$, or for every bracket $\{a'_j,b'_j\}$  in $B\cap V$, $j\not=m$, $J_k=J_{k-1}$ and $G_{k}=G_{k-1}\ast E_{p,q} $, where $p=a'_m-n$ and $q=c'_j-n$ or $q=a'_j-n$. When all blocks of $R$ have been considered, we get by construction\[J_t=s\ast H_1\ast\cdots\ast H_{2k-1}\ast s'',\quad G_t=r\ast T_1\ast\cdots\ast T_{2k-1}\ast r',\]where $T_{2i-1}=H_{2i-1}$ of $T_{2i-1}=E_{2i-1}$, and $s''$ is the resulting product of the permutations $s_i$, while $r\in\R\Sym_n$ is the product of the firstly defined $E_{p,q}$ by $s$, and $r'\in \R\Sym_n$ is the product of $s''$ by the lastly considered $E_{p,q}$ in each block.
\end{proof}

\subsection{\it Proof of Theorem \ref{PresenPBrn}}

The goal in this subsection is to prove  Theorem \ref{PresenPBrn}. 

For a monoid $M$ and a generating set $X$ of it, we denote by $\equiv_M$ the congruence on the free monoid $X^{\ast}$ generated by the defining relations of $M$ respect to $X$, i.e. $M=X^{\ast}/\equiv_M$.

In what follows, $X_n$ denotes the set of defining generators of $\Q_n$.

\begin{lemma}\label{epi}
The mapping $\S_i\mapsto\widetilde{L}_i$, $\t_i\mapsto\widetilde{H}_i$, $\F_i\mapsto\widetilde{F}_i$ and $\E_i\mapsto\widetilde{E}_i$ from $X_n$ to $\R\Br_n$ induces a monoid epimorphisms $\psi$ from  $X_n^{\ast}$ to $\R\Br_n$ and another one $\varphi$ from $\Q_n$ to $\R\Br_n$.
\end{lemma}
\begin{proof}
Lemma \ref{GenPBrn} says that $\psi$ is an epimorphism. A direct checking shows that the generators $\widetilde{L}_i$, $\widetilde{H}_i$, $\widetilde{E}_i$, $\widetilde{F}_i$ satisfy the defining relations of $\Q_n$. Then the proof that $\varphi$ is an epimorphism follows from Lemma \ref{GenPBrn}.
\end{proof}

\begin{definition}
Given a word $u\in X_n^{\ast}$, we denote by $\overline{u}$ the word obtained from $u$ by setting $\E_i=1$ and $\F_i=\t_i$, for all generators $\E_i$ and $\F_i$ occurring in $u$.
\end{definition}

\begin{remark}\label{IR}
Having in mind Remark \ref{RamDoub}, we deduce that if $\psi(u)=(I,R)$ for some $u\in X_n^{\ast}$, then $\psi(\overline{u})=(I,I)$.
\end{remark}

\begin{definition}[Diagrams]
Set the {\it diagrams} of $\S_i$ and $\t_j$ in $\Q_n$ to be the diagrams of $L_i$ and $H_j$ respectively, as elements of the partition monoid, see Figure \ref{Fig1}. The {\it diagram} of $\E_i$ in $\Q_n$ will be the diagram of it as an element of $\P_n$, see Figure \ref{003}. The {\it diagram} of $\F_i$ in $\Q_n$ is obtained by adding a tie connecting the brackets of $H_i$, see Figure \ref{FigRel1}. More generally, the {\it diagram} $D(u)$ of a word $u\in X_n^{\ast}$ is obtained by taking its generators, one by one from left to right, and connecting their diagrammatic representations from top to bottom.
\begin{figure}[H]
\[\vcenter{\hbox{
\includegraphics{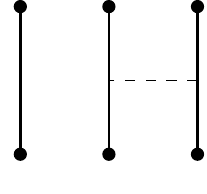}
}}
\qquad\qquad
\vcenter{\hbox{
\includegraphics{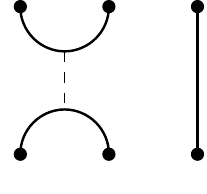}
}}
\]
\caption{For $n=3$, the diagrams of the generators $\E_2$, $\F_1$.}
\label{FigRel1}
\end{figure}
\end{definition}

\begin{remark}
For every word $u\in X_n^{\ast}$, $D(u)$ contains $n$ arcs and possibly some closed arcs. Observe that every tie in $u$ coming from a generator $\E_i$ or from a generator $\F_i$, connects two pieces of different arcs or of the same arc. Further, note that the non-tie arcs behave just like in the partition monoid, so we just need to characterize how the ties are connecting those arcs.
\end{remark} 

We define an equivalence in the set of arcs $D(u)$, including the closed arcs.

\begin{definition}
For $u\in X_n^{\ast}$, two arcs $a$ and $a'$ of $D(u)$ are said {\it tie--connected}, if a tie connects one to the other, or if there  are $k>2$ arcs $a_1,\ldots,a_k$ such that $a_1=a$, $a_k=a'$ and $a_i$ to $a_{i+1}$ are tie--connected, for all $i=1,\ldots,k-1$.
\end{definition}

\begin{proposition}\label{classes}
Let $D$ be the diagram of $u\in X_n^{\ast}$ and put $\psi(u)=(I,R)$. Two arcs of $D$ belong to the same tie-class if and  only if the  elements of $[n]\cup[n']$ joined by the corresponding arcs of $I$ belong to the same block of $R$.
\end{proposition}
\begin{proof} 
After having verified that the statement is true for the diagrams corresponding to the four types of generators,  suppose that the statement is true for any  word $u$  of length $m$ with $\psi(u)=(I,R)$. We will prove the statement for the element $ug$, where $g\in X_n$. By the induction hypothesis, a pair $j',k'$ belongs to the same block of $R\cap[n']$ if and only if $j'$ and $k'$ are endpoints of arcs tie--connected. Let $D'$ be the diagram of $ ug$ and $\psi(ug)=( I',R')$.

\begin{enumerate}
\item[$g=\E_i$.]The endpoints of the arcs of $D'$ are the same as those of $D$, since $\psi(\E_i)=(1,E_i)$. Therefore, it is evident that the tie--classes of $D'$ are still in bijection with the   blocks of $R'=R\ast E_i$.
\item[$g=\S_i$.]The endpoints of the arcs of $D'$ are the same of those of $D$, but the endpoints $i'$ and $i'+1$  which are exchanged. Since $\psi(\S_i)=(L_i,L_i)$, also $R'=R\ast L_i$ takes into account such a transposition. Therefore, the two arcs in $D'$ still are tie--connected if and only if they belong to the same  block of $R'$.
\item[$g=\t_i$.]The endpoints of the arcs of $D'$ are the same of those of $D$, but the endpoints $i'$ and $i'+1$. Indeed, these two endpoints of $D$ are connected by an arc in $D'$ (the top arc of the tangle $\t_i$), while the new endpoints $i'$ and $i'+1$ of $D'$ become the endpoints of a new arc. This is reflected by the fact that $ I'= I \ast H_i$. Therefore, two arcs of $D'$ are tie--connected if and only if they belong to the same block of $R'= R\ast H_i$.
\item[$g=\F_i$.]As for the arcs of $D$ and $D'$, the situation is as the preceding one since $I'= I \ast H_i$. However, in this case the new born arc of $D'$ with endpoints $i'$ and $i'+1$, is tie--connected to the arc obtained by connecting the endpoints $i'$ and $i'+1$ of $D$ by means of the top arc of the tangle $\F_i$. Therefore, two arcs of $D'$ are tie--connected if and only if they belong to the same block of $R'=R\ast E_i$.
\end{enumerate}
\end{proof}

\begin{definition}\label{qe}
Given $u\in X_n^{\ast}$, we define the word $u^{\E}\in X_n^{\ast}$ obtained from $u$ through the following steps:
\begin{enumerate}
\item[Step 1.]Replace every $\t_i$ by $\E_i\t_i\E_i$, and every $\F_i$ by $\E_i\F_i\E_i$. See (\ref{Eiti}) and (\ref{EiFi}).
\item[Step 2.]Double every tie by  means of relations (\ref{Ei2}) and (\ref{Fi2}) and move it at  the left (if possible) by using (iii)--(iv) Lemma \ref{ties}; the same to move at the right. Repeat this step for every extended tie obtained, until no new ties can be added.
\item[Step 3.]Finally, replace by $\F_k$ every $\t_k$ such that it is preceded and followed by $\E_{i,k}$ for some $i$ or such that it is preceded and followed by $\E_{k+1,j}$ for some $j$.
\item[Step 4.]For each $\F_k$ considered at step 3, apply step 2, point (iv) Lemma \ref{ties} in particular, to those ties preceding and following it, that before step 3 have stopped to move because of point (v) Lemma \ref{ties}.
\end{enumerate}
\end{definition}
\begin{figure}[H]
\[\begin{array}{c}
\vcenter{\hbox{
\includegraphics{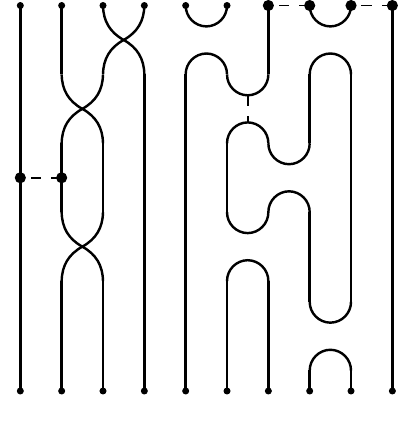}
}}
\stackrel{\text{step 1}}{\longrightarrow}
\vcenter{\hbox{
\includegraphics{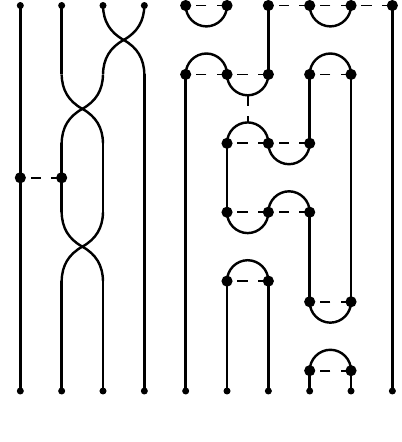}
}}
\stackrel{\text{step 2}}{\longrightarrow}
\\[2.5cm]
\vcenter{\hbox{
\includegraphics{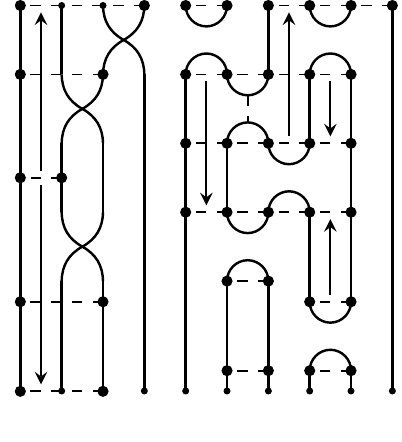}
}}
\stackrel{\text{step 3}}{\longrightarrow}
\vcenter{\hbox{
\includegraphics{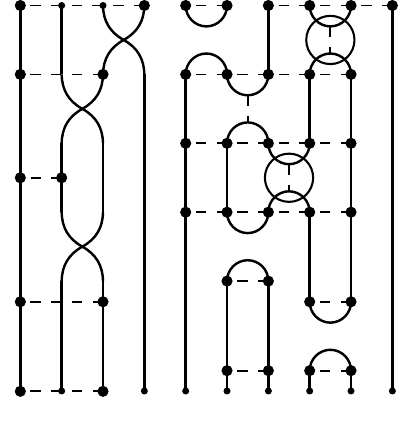}
}}
\stackrel{\text{step 4}}{\longrightarrow}
\vcenter{\hbox{
\includegraphics{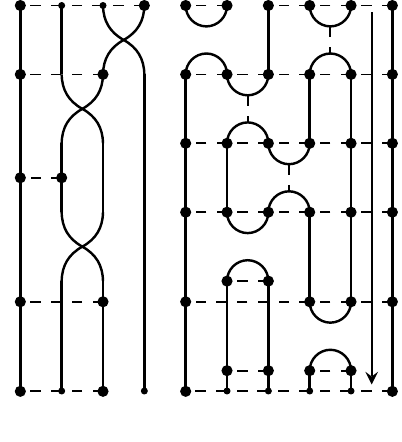}
}}
\end{array}\]
\caption{For $n=10$, the construction of $(\S_3\t_5\t_8\S_2\F_6\E_1\t_7\S_2\t_6)^{\E}$.}
\label{d20}
\end{figure}
 
We denote by $\sim $ the equivalence relation on $X_n^{\ast}$ generated by (\ref{P2}) and (\ref{P3}).

\begin{proposition}\label{qqe}
Let $u,v\in X_n^{\ast}$, we have: $u^{\E}\sim v^{\E}$ if and only if the tie--classes of $D(u)$ and $D(v)$ coincide.
\end{proposition}
\begin{proof}
Firstly, observe that every tie $\E_{i,j}$ added to $u$ to form $u^\E$, lies between arcs already tie--connected. The same for every tie added by replacing $\t_i$ by $\F_i$. Therefore we can extend the definition of arcs tie--connected to the extended ties without affecting the tie--classes of the diagram $D(u)$, that will coincide with those of $D(u^{\E})$. Secondly, notice that, while a tie $\E_i$  occurring in $u$, by the procedure explained in Definition \ref{qe}, moves at the left and right in $u$, the corresponding tie in $D(u)$ moves up and down. If the original tie lies between two arcs that are lines in $\psi(u)$, then the copies of it reached the top and the bottom of the diagram. Otherwise, the tie may stop to move only in presence of an element $\t_j$, with which it does not commute. Moreover, by construction, all pieces of the same arcs result to be tied. Finally, observe that any closed curve appearing in $D( u)$, is filled by ties, and every tangle belonging to it is a tied tangle. 
\end{proof}

\begin{remark}\label{e1}
Note that by setting $\E_i=1$, we get $E_{i,j}=1$ for every $j$, see Definition \ref{extended}. 
 Thus $\overline{u}=\overline{u^{\E}}$.
\end{remark}
  
\begin{lemma}\label{overline}
Let $u,v\in X_n^{\ast}$ satisfying $\overline{u}\equiv_{\Q_n}\overline{v}$. Then $u\equiv_{\Q_n}v$ if and only if $\psi(u)=\psi(v).$
\end{lemma}
\begin{proof}
By Lemma \ref{epi}, if $u\equiv_{\Q_n}v$ then $\psi(u)=\psi(v)$. Reciprocally, suppose that $\psi(u)=\psi(v)$. By Proposition \ref{classes} and Proposition \ref{qqe}, we have $ u^\E\sim v^{\E}$. Hence $u\equiv_{\Q_n}v$.
\end{proof}

\begin{lemma}\label{X}
Let $u,v\in X_n^{\ast}$ with no letters $e_i,f_j$ occurring on them, such that $u\equiv_{\Br_n}v$, and let $I\in\mathfrak{B}_n$ such that $\psi(u)=\psi(v)=(I,I)$. Then, there are exactly $b_n$ distinct relations $\mathcal{R}^j:p_j\equiv_{\Q_n}q_j$, where $\overline{p}_j\equiv_{\Br_n}u$, $\overline{q}_j\equiv_{\Br_n}v$ and $\psi(p_j)=\psi(q_j)=(I,R_j)$ realizes every ramified partition of type $(I,R)$.

\end{lemma}
\begin{proof}
The existence of the epimorphism $\psi$, guarantees that the number of distinct relations $\mathcal{R}^j$ is at least $b_n$. However, if this number exceeds $b_n$, there is a pair $w,w'\in X_n^*$ satisfying $w\equiv_{\Q_n}w'$ such that $\overline{w}\equiv_{\Br_n}u$ and $\psi(w)=\psi(p_j)$ for some $j$ but $w\not=p_j$ for all $j$. This contradicts Lemma \ref{overline}.
\end{proof}

In fact, every defining relation $u\equiv_{\Br_n}v$ of $\Br_n$ involves at most two different indices of generators, and consequently a number $k\leq4$ of arcs of $\psi(u)=(I,I)$. For the remaining $n-k$ arcs, $I$ coincides with the identity. Thus, neglecting the vertical arcs, there are exactly $b_k$ distinct relations $\mathcal{R}^j$ in $\Q_n$ obtained from $u$ and $v$  by changing some $\t_i$ into $\F_i$, and by inserting some tie $\E_{i,j}$ between the $k$ arcs.

\begin{example}
Consider the relation $\S_1\t_2\S_1=\S_2\t_1\S_2$ in $\Br_n$. It involves two indices, $k=3$ and $b_k=5$. The five relations $\mathcal{R}^j:p_j=q_j$ in $\Q_n$ corresponding to the double partitions $(I,R_j)$, where $ I =\{\{1,3\},\{2,2'\},\{1',3'\}\}$ are:\[\begin{array}{lrcll}
\mathcal{R}^1:&\E_{1,2}\S_1\E_2\t_2\E_2\S_1\E_{1,2}&=&\E_{1,2}\S_2\E_1\t_1\E_1\S_2\E_{1,2},
&R_1=\{\{1,3\},\{2,2'\},\{1',3'\}\},\\
\mathcal{R}^2:&\E_1\E_2\S_1\E_1\E_2\t_2\E_2\S_1\E_{1,2}&=&\E_1\E_2\S_2\E_1\E_2\t_1\E_1\S_2\E_{1,2},
&R_2=\{\{1,2,3,2'\},\{1',3'\}\},\\
\mathcal{R}^3:&\E_{1,2}\S_1\E_2\t_2\E_1\E_2\S_1\E_1\E_2&=&\E_{1,2}\S_2\E_1\t_1\E_1\E_2\S_2\E_1\E_2,
&R_3=\{\{1,3\},\{2,2',1',3'\}\},\\
\mathcal{R}^4:&\E_{1,2}\S_1\E_2\F_2\E_2\S_1\E_{1,2}&=&\E_{1,2}\S_2\E_1\F_1\E_1\S_2\E_{1,2},
&R_4=\{\{1,3, 1',3'\},\{2,2'\}\},\\
\mathcal{R}^5:&\E_1\E_2\S_1\E_1\E_2\F_2\E_1\E_2\S_1\E_1\E_2&=&\E_1\E_2\S_2\E_1\E_2\F_1\E_1\E_2\S_2\E_1\E_2,
&R_5=\{\{1,2,3,1',2',3'\}\}.
\end{array}\]

\begin{figure}[H]
\[\begin{array}{c}
\vcenter{\hbox{
\includegraphics{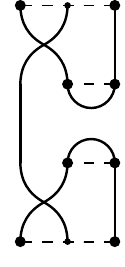}
}}=\!\vcenter{\hbox{
\includegraphics{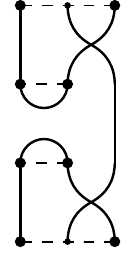}
}}\qquad\vcenter{\hbox{
\includegraphics{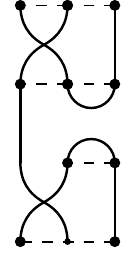}
}}=\!\vcenter{\hbox{
\includegraphics{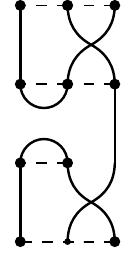}
}}\qquad\vcenter{\hbox{
\includegraphics{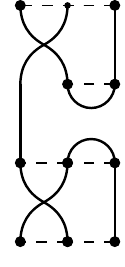}
}}=\!\!\vcenter{\hbox{
\includegraphics{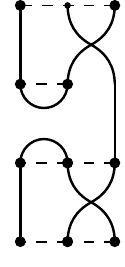}
}}\\[2cm]
\vcenter{\hbox{
\includegraphics{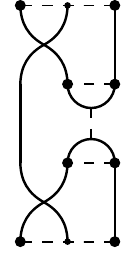}
}}=\!\vcenter{\hbox{
\includegraphics{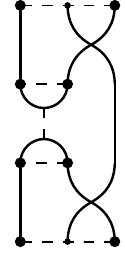}
}}\qquad\vcenter{\hbox{
\includegraphics{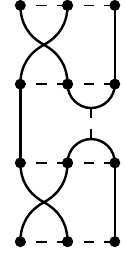}
}}=\!\vcenter{\hbox{
\includegraphics{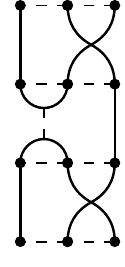}
}}
\end{array}\]
\caption{The five relations corresponding to $\S_1\t_2\S_1=\S_2\t_1\S_2$.}
\label{d15}
\end{figure}
\end{example}

\begin{lemma}\label{inj} 
Let $u,v\in X_n^{\ast}$ satisfying $\psi(u)=\psi(v)$. Then $u\equiv_{\Q_n}v$.
\end{lemma}
\begin{proof}
We will prove that $u^{\E}\equiv_{\Q_n}v^{\E}$. Remark \ref{IR} implies that $\psi(\overline{u})=\psi(\overline{v})$, so by Lemma \ref{BLBr} we get $\overline{u}\equiv_{\Br_n}\overline{v}$. Therefore, there is a sequence of defining relations $R_i:r_i\equiv_{\Br_n}r_{i+1}$ with $1\leq i\leq t$ such that $\overline{u}=w_1\equiv_{\Br_n} w_2\equiv_{\Br_n}\cdots\equiv_{\Br_n}w_t=\overline{v}$. We will define a corresponding sequence of relations in $\Q_n$, transforming $u^{\E}$ into $v^{\E}$. For $i=1$ we have $\overline{u}=w_1=ar_1b$ and $w_2=ar_2b$ for some words $a,b$. Observe that $u^{\E}=\hat{w}_1\equiv_{\Q_n}\hat{w}_2$, where $\hat{w}_1=ct_1d$ and $\hat{w}_2=ct_2d$, in such a way that $\overline{t}_1=r_1$, $\overline{t}_2=r_2$, $c$ is a word that ends with a generator different from a tie and $d$ is a word that begins with a generator different from a tie. The product of $\E_{i,j}$ at the beginning of $t_1$ allows us to choose the relation $R_1^j$ by using Lemma \ref{X}, and replace $t_1$ with the corresponding $t_2$. Observe that this replacement preserves the products of $\E_{i,j}$ at the beginning and at the end of $t_1$. We proceed in this way until we get $u^{\E}=\hat{w}_1\equiv_{\Q_n} \hat{w}_2\equiv_{\Q_n}\cdots\equiv_{\Q_n}\hat{w}_t=v^{\E}$.
\end{proof}

\begin{theorem}\label{PresenPBrn}
The epimorphism $\varphi$ of Lemma \ref{epi} is an isomorphism.
\end{theorem}
\begin{proof}
It follows from Lemma \ref{epi} and Lemma \ref{inj}.
\end{proof}

\begin{corollary}\label{normalformQ}
Every element $q\in\Q_n$ either belongs to $\T\Sym_n$ or for some $k\le n/2$ can be  written in the following normal form\[q_N=r\ v_1v_3\cdots v_{2k-1}\ r',\]where $r,r'\in\T\Sym_n$ and $v_{2i-1}=\t_i$ or $v_{2i-1}=\F_i$.
\end{corollary}
\begin{proof}
It follows from Theorem \ref{PresenPBrn} and Lemma \ref{GenPBrn}.
\end{proof}

\section{The monoid $\bBr_n$}\label{sectionsix}

We give a presentation of the submonoid $\bBr_n$ of $\R\Br_n$ formed by those ramified partitions $(I, R)$ of $\R\B_n$ with $R$ {\it balanced}, see Definition \ref{BalPart}. As we did with $\R\Br_n$, this presentation is obtained by proving that $\bBr_n$ is isomorphic to an auxiliary monoid $W_n$, which is defined from the presentation of $\Q_n$ by omitting in the the generators $\t_i$'s together with relations containing them, and adding certain relations valid in $\Q_n$. Also, the cardinality of $\bBr_n$ is calculated.

\begin{definition}\label{BalPart}
In a ramified partition $(I,R)$ of $\R\Br_n$, a block of the partition $R$ is said balanced if it contains the same number (possibly zero) of up brackets and down brackets from $I$. We say that $R$ is balanced if each block of it is balanced.
\end{definition}

We will denote by $\bBr_n$ the set of ramified partitions $(I,R)\in\R\Br_n$ such that $R$ is balanced.

Let $\W_n$ be the monoid presented by generators $\S_1,\ldots,\S_{n-1}$, $\E_1,\ldots,\E_{n-1}$, $\F_1,\ldots,\F_{n-1}$, satisfying the relations (\ref{S1}) to (\ref{S3}), (\ref{Fi2}) to (\ref{FiFjFi}) together with the following relations:
\begin{align}
\S_i\E_i&=\E_i\S_i\quad\text{for all }i,\label{tSiEi}\\
\S_i\E_j&=\E_j\S_i\quad\text{for all $i$ and $j$ s.t. }\vert i-j\vert>1,\label{tSiEj}\\
\E_i\E_j\S_i&=\S_i\E_i\E_j=\E_j\S_i\E_j\quad\text{for all $i$ and $j$ s.t. }\vert i-\vert=1,\label{tEiEjSi}\\
\E_i\S_j\S_i&=\S_j\S_i\E_j\quad\text{for all $i$ and $j$ s.t. }\vert i-j\vert=1,\label{tEiSjSi}\\
\S_i\F_i&=\F_i\S_i=\F_i\quad\text{for all }i,\label{tSiFi}\\
\S_i\F_j&=\F_j\S_i\quad\text{for all }\vert i-j\vert>1,\label{tSiFj}\\
\F_i\F_j\S_i&=\E_j\F_i\S_j\quad\text{for all $i$ and $j$ s.t. }\vert i-j\vert=1,\label{tFiFjSj}\\
\S_j\F_i\F_j&=\S_i\F_j\E_i\quad\text{for all $i$ and $j$ s.t. }\vert i-j \vert=1,\label{tSjFiFj}\\
\F_i\S_j\F_i&=\E_j\F_i\E_j\quad\text{for all $i$ and $j$ s.t. }\vert i-j\vert=1,\label{tFiSjFicasetBrn}\\
\S_i\F_j\S_i&=\S_j\F_i\S_j\quad\text{for all $i$ and $j$ s.t. }\vert i-j\vert=1,\label{tSiFjSi}\\
\S_i\S_j\F_i&=\F_j\S_i\S_j\quad\text{for all $i$ and $j$ s.t. }\vert i-j\vert=1,\label{tSiSjFi} \\
\E_i\S_j\F_i&=\S_j\F_i\E_j\quad\text{for all $i$ and $j$ s.t. }\vert i-j\vert=1,\label{tEiSjFi}\\
\F_i\S_j\E_i&=\E_j\F_i\S_j\quad\text{for all $i$ and $j$ s.t. }\vert i-j\vert=1.\label{tFiSjEi}
\end{align}

\begin{proposition}\label{submonoid}
The monoid $\W_n$ is a submonoid of $\Q_n$.
\end{proposition}
\begin{proof}
The monoid $\W_n$ is generated by the same generators as $\Q_n$ excluding the tangle generators $\t_i$. All the defining relations of $\W_n$ hold in $\Q_n$, and every relations of $\Q_n$ non involving the generators $\t_i$'s holds in $\W_n$.
\end{proof}

To show that $\varphi(\W_n)=\bBr_n$, we will use the following lemma.

\begin{lemma}\label{psiW}
If $w\in\W_n$ such that $\varphi(w)=(I,R)$, then $R$ is balanced.
\end{lemma}
\begin{proof}
The statement is true when $w$ is one of the generators. Suppose that it holds for any word $w$ being the product of $m-1$ generators. We prove that it holds also for the word $wg$ with $m$ generators. Let $\varphi(w)=(I,R)$ and $\varphi(wg)=(I',R')$. If $g=\t_i$, the numbers of up brackets and down brackets in the blocks of $R$ and of $R'$ are the same. If $g=\E_i$, $I'=I$ and there are two possibilities: if $i'$ and $i'+1$ belong to the same block of $R$ or not. In the first case we have, then $R'=R$. In the other case, these blocks merge in a unique block of $R'$. In it, the number of up brackets is equal to the number of down brackets, since by hypothesis these numbers coincide in the merging blocks. Finally, if $g=\F_i$, there are four possibilities: (1) if $i'$ and $i'+1$ belong to two different lines $\{k,i'\}$ and $\{j,i'+1\}$ of $I$, in $I'$ we get the new up bracket $\{k,j\}$ and the new down bracket $\{i',i'+1\}$. These two new brackets lie in the same block of $R'$, that by consequence results to be balanced. (2) If $i'$ and $i'+1$ belong to a down bracket in $I$, then in $R$ there is a block containing it together with $r\ge0$ down brackets and $r+1$ up brackets by hypothesis. In $I'$ the down bracket $\{i',i'+1\}$ remains in the same block as in $R$, still balancing the number of up brackets in it. (3) If $i'$ belongs to a line $\{k,i'\}$ and $i'+1$ to a down brackets $\{i'+1,j'\}$ of $I$, then we get in $I'$ a line $\{k,j'\}$ and a down bracket $\{i',i'+1\}$. Such a bracket is in the same block of $R'$ as the line $\{k,j'\}$. This line shares the block with an up bracket of $I'$, already existing in $I$, that balanced in a block of $R$ the sum of up and down brackets together with the disappeared brackets $\{i'+1,j'\}$. Therefore $R'$ is still balanced. The case in which $i'$ belongs to a down brackets and $i'+1$ to a line, is similar to this last case. (4) If $i'$ and $i'+1$ belong to two down brackets of $I$, i.e. $b_1:=\{i', j'\}$ and $b_2:=\{i'+1,k'\}$, we get in $I'$ two down brackets $b_3:=\{j',k'\}$ and $b_4:=\{i',i'+1\}$. Suppose that, in $R$, $b_1$ and $b_2$ belong to two blocks $B_1$ and $B_2$, not necessarily distinct. Then we get in $R'$ a unique block $\B':=\{B_1\cup B_2\cup b_3\cup b_4\}\setminus\{b_1,b_2\}$, which is balanced since $B_1$ and $B_2$ are balanced by hypothesis.
\end{proof}

\begin{theorem}\label{isoW}
We have $\varphi(\W_n)=\bBr_n$. Thus, the definition of $\W_n$ is a presentation for $\bBr_n$.
\end{theorem}
\begin{proof}
By Lemma \ref{psiW}, we get that $\varphi(\W_n)\subseteq\bBr_n$. So, the proof is finished by proving that for every $(I,R)\in\bBr_n$ there is $w\in\W_n$ such that $\varphi(w)=(I,R)$. More precisely, given $(I,R)\in\bBr_n$, we show how to write it in terms of the $\widetilde{L}_i$'s, $\widetilde E_i$'s and $\widetilde F_i$'s. Then the proof ends by substituting the generators one by one by their preimages under $\varphi$, i.e., the $\S_i$'s, $\E_i$'s and $\F_i$'s.
 
We order the blocks $B_j$ of $R$ by the minimal left endpoint of the up brackets in them. The blocks with no brackets are put at the end, and are ordered by the minimal upper endpoint of the lines in it. Then, for each block $B_j$ we proceed this way. Let $\{a_1^{(j)},b_1^{(j)}\},\dots,\{a_{k_j}^{(j)},b_{k_j}^{(j)}\}$ be the up brackets and $\{{a'_1}^{(j)},{b'_1}^{(j)}\},\dots,\{{a'_{k_j}}^{(j)},{b'_{k_j}}^{(j)}\}$ the down brackets in $B_j$ satisfying $a_i<b_i$, $a_{i}<b_{i+1}$, $a'_i<b'_i$, $a'_{i}<b'_{i+1}$. Moreover, let $\{c_1^{(j)},{c'_1}^{(j)}\},\dots,\{c_{l_j}^{(j)},{c'_{l_j}}^{(j)}\}$ be the lines of $B_j$, satisfying $c_i<c_{i+1}$. Let $m_j=2k_j+l_j$ be the cardinality of the block $B_j$. We define the permutations:
\[s:=\left(\begin{array}{cccccccccccc}
a_1^{(1)}&b_1^{(1)}&\dots&a_{k_1}^{(1)}&b_{k_1}^{(1)}&c_1^{(1)}&\dots&c_{l_1}^{(1)}&a_1^{(2)}&b_1^{(2)}&\dots&c_{l_t}^{(t)}\\
1&2&\dots&2k_1-1&2k_1&2k_1+1&\dots&m_1&m_1+1&m_1+2&\dots&n\end{array}\right),\]
\[s':=\left(\begin{array}{cccccccccccc}
1&2&\dots&2k_1-1&2k_1&2k_1+1&\dots&m_1&m_1+1&m_1+2&\dots&n\\
{a'}_1^{ (1)}&{b'}_1^{(1)}&\dots&{a'}_{k_1}^{(1)}&{b'}_{k_1}^{(1)}&{c'}_1^{(1)}&\dots&{c'}_{l_1}^{(1)}&{a'}_1^{(2)}&{b'}_1^{(2)}&\dots&{c'}_{l_t}^{(t)}
\end{array}\right).\]
Also, set $M_1=0$ and $M_j=M_{j-1}+m_{j-1}$ for $j=2,\dots,t$, and define:\[\mathbb{E}=\prod_{j=1}^tE^{(j)},\quad E^{(j)}:=\prod_{i=1}^{m_j-1}\widetilde{E}_{i+M_j },\]\[\mathbb{F}=\prod_{j=1}^tF^{(j)},\quad F^{(j)}:=\prod_{i=1}^{k_j}\widetilde{F}_{2i-1+M_j}.\]Using the same arguments as in Lemma \ref{GenPBrn}, we obtain\begin{equation}\label{normalW}(I,R)=s\mathbb{E}\mathbb{F}s',\end{equation}where $s_1$ and $s'_1$ are written in terms of the $\widetilde L_i$'s.
\end{proof} 

Finally, we are going to calculate the cardinality of $\bBr_n$. A 2-balanced partition is a partition of a set of in which there are elements of three types, say positive, negative and neutral, with the condition that in  each block the number of positive elements equals the number of the negative elements. Evidently, the set contains the same number, say $k$, of positive and negative elements. We denote $U(n,k)$ the number of 2-balanced partitions of a set with $n$ elements, whose $k$ are positive and $k$ are negative. Therefore, $U(n,k)$ is the number of balanced partitions $R$ such that the number of up and down arcs of $I$ are exactly $k$. Evidently, $U(n,0)=b_n$.

The numbers $U(n,k)$ are given by the triangle in OEIS A343254.

\begin{proposition}
The cardinality of $\bBr_n$ is\begin{equation}\label{dimens}\sum_{k=0}^{n/2}\frac{n!^2 }{2^{2k}k!^2(n-2 k)!}\ U(n,k).\end{equation}
\end{proposition}
\begin{proof}
Fix the number $k$ of pairs of up and down arcs. The number of partitions $I$ with such value of $k$ is obtained by multiplying the number of choices of the endpoints of the up arcs, $\frac{n!}{2^k k!(n-2k)!}$, by the same number of choices of the endpoints of the down arcs, and by the number of choices, $(n-2k)!$, of the lower points for the remaining $n-2k$ lines.
\end{proof}

Here the list of sizes of $\bBr_n$ for $n=1,\dots,14$:
\[\begin{tabular}{|c|c|}\hline
1&1\\\hline
2&5\\\hline
3&48\\\hline
4&747\\\hline
5&17040\\\hline
6&531810\\\hline
7&21634515\\\hline
\end{tabular}\qquad
\begin{tabular}{|c|c|}\hline
8&1107593235\\\hline
9&69482175840\\\hline
10&5229801016650\\\hline
11&464302838867175\\\hline
12&47939037056237250\\\hline
13&5688447976340254125\\\hline
14&767923605609975114300\\\hline
\end{tabular}\]

\section{The monoid $\tJ_n$ and boxed partitions}\label{sectionsev}

Here, we consider a submonoid $\tJ_n$ of $\Q_n$ generated by the $\E_i$'s and the $\F_i$'s subject to certain relations that resemble the defining relations of $\J_n$, in such a way that it can be understood as a tied Jones monoid. We show a normal form for the elements of $\tJ_n$. Further, we show that $\tJ_n$ can be realized as the submonoid of $\bBr_n$ formed by $(I,R)$'s, where $R$ is a boxed set partition (Definition \ref{DefBoxpart}).

\subsection{}

We define $\tJ_n$ as the monoid presented by generators $\E_1,\ldots,\E_{n-1}$, $\F_1,\ldots,\F_{n-1}$ subject to the relations (\ref{Fi2})--(\ref{FiFjFi}).

Observe that by setting $\E_i=1$, for all $i$, the monoid $\tJ_n$ becomes the Jones monoid $\J_n$. In other words, the mapping $\F_i\mapsto\t_i$, $\E_i\mapsto1$ defines a monoid homomorphism from $\tJ_n$ to $\J_n$. This remark together with the fact that the $\E_i$'s are tie elements, and the $\F_i$'s are tied tangles, allows to say that $\tJ_n$ is a tied version of the Jones monoid. See \cite{AiJuPa2019} for other tied versions of $\J_n$. As for the Jones monoid, it is shown that the elements of $\tJ_n $ have a normal form. To do this, we need to introduce the following elements:\begin{equation}\label{Frs}
\F_{r,s}:=\F_s\F_{s-1}\cdots\F_r,\quad\text{where }0<r\leq s< n.
\end{equation}
Observe that $\F_{r,r}:=\F_r$.

\begin{proposition}[Normal form]\label{proEF}
Every element of $\tJ_n $ can be written uniquely in the form $\F\E$, with
\begin{equation}\label{EF}
\F:=\F_{j_1,k_1}\cdots\F_{j_{t},k_{t}},\quad\quad\E:=\prod_k\E_{i_k},
\end{equation}where $0<j_i<j_{i+1}<n, 0<k_i<k_{i+1}<n$ and each index $i_k$ is absent in $\F$.
\end{proposition}
\begin{proof}
Relations (\ref{EiEj}) and (\ref{EiFj}) imply that the tie generators commute with all other generators; therefore we can put every $\E_i$ to the right of the word. So, every $w$ in $\tJ_n$ can be written in the form $w=\F\E$, where $\F$ is a word in the $\F_i$'s, and $\E$ is a word in the $\E_i$'s,  which is taken in normal form according to Proposition \ref{FNPn}. We can proceed now as Jones did in \cite[Lemma 4.1.2]{JoIM1983}, but on the word $\F$, to reduce it to a product as in (\ref{EF}). Observe that in our situation the reduction of $\F_i\F_j\F_i$ involves generators $\E_j$, see relation (\ref{FiFjFi}). However, these tie generators can always be moved to the right. Therefore the proof follows arguing as in \cite[Aside 4.1.4]{JoIM1983}.

Finally, thanks to the relation (\ref{EiFi}), every generator $\E_j$ can be removed from $\E$ if $\F_j$ occurs  in $\F$.
\end{proof}

As an application of Proposition \ref{proEF} we will compute the size of $\tJ_n$. For this porpuse we  need to introduce below some notations.

Let $\F$ be a product of $\F_{j_i,k_i}$'s as in (\ref{EF}). We define:
\begin{enumerate}
\item$N(\F)$ as the number of different indices of the generators $\F_i$ in $\F$.
\item$\mathcal{G}_n^k$ as the set formed by the elements $\F$ with $N(\F)=k$. We denote by $G_n^k$ the cardinality of $\mathcal{G}_n^k$.
\item We call {\it gap} an index $g$ that in expression (\ref{EF}) satisfies $k_i<g<j_{i+1}$ for some $i$ such that $0<i<t$.
\end{enumerate}

\begin{example}
For $\F=\F_{1,2}\F_{4,5}=\F_2\F_1\F_5\F_4\in\tJ_6$, we have $N(\F)=4$ and  that $g=3$ is a gap.
\end{example}

The next lemma says that the integers $G_n^k$ form the sequence known as the Catalan triangle $T(n,k)$, that is, the sequence defined by initial conditions: $T(0,0)=1$, $T(n,0)=1$, and $T(n,n)=0$ for $n>0$, together with the following recursive formula\[T(n,k)=T(n,k-1)+T(n-1,k),\quad\text{for }n>1,\,0<k<n.\]
For more details see \cite{RoDM1978}.

\begin{lemma}\label{GCatalan}
For every $n\geq0$ and $k\leq n$, we have $G_n^k=T(n,k)$.
\end{lemma}
\begin{proof}
Note that $N(\mathbf{1})=0$, so $\mathcal{G}_n^0=\{\mathbf{1}\}$ and $G_n^0=1$. Moreover, $G^n_n=0$ for $n>0$ since $\tJ_n$ has $n-1$ different generators $\F_i$. Thus, the proof is concluded if we prove that $G_n^k$ satisfies the recursive formula of the Catalan triangle, i.e.
\begin{equation}\label{recurr}
G_n^k=G_{n-1}^k+G_{n}^{k-1}.
\end{equation}
It is evident that $\mathcal{G}_{n-1}^k\subset\mathcal{G}_n^k$. Therefore we have to prove that there is a  bijection $h$ between $\mathcal{G}_{n}^k\setminus\mathcal{G}_{n-1}^k$ and $\mathcal{G}_{n}^{k-1}$. Observe that an element in $\mathcal{G}_{n}^k\setminus\mathcal{G}_{n-1}^k$ contains $\F_{n-1}$, and contains it only once, while  an element in $\mathcal{G}_{n}^{k-1}$ may contain or not $\F_{n-1}$. Consider an element $\Bf\in\mathcal{G}_{n}^{k-1}$. If it does not contain $\F_{n-1}$, then define $h(\Bf)=\Bf\F_{n-1}$. Otherwise, $\Bf$ terminates in $\F_{j_t,n-1}$. Observe that $t\le k$. Let $g$ be the maximum gap in $\Bf$. Observe that $g<n-1$ and, by definition, there is an index $j_i$ such that $k_{i-1}<g=j_{i}-1$. We define in this case\[h(\Bf)=\F_{j_1,k_1}\cdots\F_{j_{i-1},k_{i-1}}\ \F_{g,k_i}\ \F_{j_{i+1}-1,k_{i+1}}\cdots\F_{j_t-1,n-1}.\]If $\Bf$ has no gaps, we set  $g=j_1-1$. 

Figure \ref{d18} shows 14 diagrams made of white boxes that represent the elements $\Bf\in\F\in\mathcal{G}_5^3$. A box at height $j$ represents $\F_j$, and each column with boxes at heights from $j$ to $k$ represents $\F_{j,k}$. An element $\Bf$ is thus represented by adjacent columns, ordered from left to right, according to the sequence of $\F_{j,k}$'s in its expression. For each diagram, the gray boxes are added to form $h(\Bf)\in\mathcal{G}^4_5$, the special index $g$ is marked.
\begin{figure}[H]
\[\begin{array}{c}
\ytableausetup{smalltableaux}
\begin{ytableau}\none&\none&\none&*(gray)\\\none&\none&&\none\\\none&&\none&\none\\&\none&\none&\none\end{ytableau}
\quad
\begin{ytableau}\none&\none&*(gray)\\\none&&\none\\&\none&\none\\&\none&\none\end{ytableau}
\quad
\begin{ytableau}\none&\none&*(gray)\\\none&&\none\\\none&&\none\\&\none&\none\end{ytableau}
\quad
\begin{ytableau}\none&\none&*(gray)\\\none&&\none\\&&\none\\&\none&\none\end{ytableau}
\quad
\begin{ytableau}\none&*(gray)\\&\none\\&\none\\&\none\end{ytableau}
\\[1.6cm]
\begin{ytableau}\none&\none&\\\none&\none&*(gray)3\\\none&&\none\\&\none&\none\end{ytableau}
\quad
\begin{ytableau}\none&\\\none&*(gray)3\\&\none\\&\none\end{ytableau}
\quad
\begin{ytableau}\none&\none&\\\none&&*(gray)\\\none&*(gray)2&\none\\&\none&\none\end{ytableau}
\quad
\begin{ytableau}\none&\none&\\\none&&*(gray)\\&*(gray)&\none\\*(gray)1&\none&\none\end{ytableau}
\quad
\begin{ytableau}\none&\\&*(gray)\\&\none\\*(gray)1&\none\end{ytableau}
\quad
\begin{ytableau}\none&\\\none&\\\none&*(gray)2\\&\none\end{ytableau}
\quad
\begin{ytableau}\none&\\\none&\\&*(gray)\\*(gray)1&\none\end{ytableau}
\quad
\begin{ytableau}\none&\\&\\&*(gray)\\*(gray)1&\none\end{ytableau}
\quad
\begin{ytableau}\\\\\\*(gray)1\end{ytableau}
\end{array}\]
\caption{Schemes representing $\Bf\in\mathcal{G}_{5}^4$. First one: $\Bf=\F_1\F_2\F_3\F_4$; last one, $\Bf=\F_{1,4}$.}
\label{d18}
\end{figure}
Clearly, we have created an element of $\mathcal{G}_n^k\setminus\mathcal{G}_{n-1}^k$ from an element of $\mathcal{G}_{n}^{k-1}$, and the definition of $h$ guarantees that $h(\Bf_1)=h(\Bf_2)$ iff $\Bf_1=\Bf_2$, so $h$ is injective. To prove the surjectivity of $h$, we consider an element $\Bf\in\mathcal{G}_{n}^k\setminus\mathcal{G}_{n-1}^k$ and we use the inverse procedure to get an element $\Bf'\in\mathcal{G}_{n}^{k-1}$. Observe that the element $\Bf$ terminates in $\F_{j_t,n-1}$. If $j_t=n-1$, then we obtain $\Bf'$ by removing   $\F_{n-1}$. Otherwise, we change $j_t$ into $j_t+1$. If $j_t+1>k_{t-1}$, then $ j_t $ is a gap, and the so obtained element is the searched $\Bf'$. Otherwise, $j_t< k_{t-1}$, therefore $j_{t-1}<k_{t-1}$ and we can change  $j_{t-1}$ into $j_{t-1}+1$. If $j_{t-1}+1>k_{t-2}$, then $j_{t-1}$ is a gap, and we have got $\Bf'$, otherwise we continue the same way till the first index $j_i$ satisfying $j_i+1>k_{i-1}$. At this point we end since we have obtained the searched $\Bf'$ with gap $g=j_i$. If there are no gaps in $\F$, we change eventually $j_1$ into $j_1+1$. Again this is possible since $j_1=k_1$ should contradicts the fact that there were no gaps before $g>j_1$.
\end{proof}

Note that in the case $k=n-1$, recurrence (\ref{recurr}) gives $G_n^{n-1}=G_n^{n-2}$, since $G_{n-1}^{n-1}=0$.

\begin{proposition}\label{dimtJ}
$\tJ_n$ has $\binom{2n-1}{n} =(2n-1)C_{n-1}$ elements.
\end{proposition}
\begin{proof}
Because of  Proposition \ref{proEF}, it is enough to count the possible normal form $ \F \E $. Now, if $ \F \in\mathcal{G}_n^k$, $\E$ is the  product of  at most $n-1-k$ different generators $\E_j$.  It follows that the total number of generators is\[\sum_{k=0}^nG_n^k2^{n-1-k}.\]This sum equals $\binom{2n-1}{n}$ by \cite[Theorem 1]{AiarXiv2020}.
\end{proof}

\subsection{}

At this point we will use the concepts of  linear partition, see Subsection \ref{subsecLP}, and the {\it boxed set partition}.

\begin{definition}\label{DefBoxpart}
A set partition $R$ of $[n]\cup[n']$ is said to be {\it boxed} if $R\cap[n]$ is linear, and for every $i\in[n]$, $i$ and $i'$ belong to a same block of $R$.
\end{definition}

\begin{remark}\label{remlinear}
Observe that an element of $\R\Br_n$, which is a product of $k$ different generators $\widetilde E_i$, defines a ramified partition $(I,R)$ where $I$ is the identity and $R$ is a boxed partition with $n-k$ blocks.
\end{remark}

\begin{definition}
We define $\bJ_n$ as the submonoid of $\bBr_n$, consisting of  the ramified partitions $(I,R)$, such that $I\in J_n$, and $R$ is boxed.

\end{definition}

We will close this section by proving the following theorem.

\begin{theorem}\label{isotj} 
The isomorphism  $\varphi$ carries $\tJ_n$ in $\bJ_n$.
\end{theorem}

In order to prove this theorem we firstly show that $\bJ_n$ and $\tJ_n$ have the same cardinality. So, let us start with some lemmas.

\begin{lemma}\label{lemmapart}
Let $\Bf=\F\E\in\tJ_n$ as in Proposition \ref{proEF} with $\psi(\Bf)=(I,R)$. If $N(\Bf)=k$ and $\E=1$, then $R$ is a boxed partition with $n-k$ blocks.
\end{lemma}
\begin{proof}
If $k=0$, $\Bf$ is the identity and $R=1$ is  boxed with $n$ blocks. Let us decompose the element $\Bf$ into the product $\Bf=\prod_q\Bf_q$, such that no one of the indices in a factor $\Bf_q$ is a gap, and there is at least a gap between the indices of two consecutive factors. I.e., each factor has the form\[\Bf_q:=\F_{j_i,k_i}\cdots\F_{j_{l},k_{l}},\]where $j_{i+1}\le k_i$. We will denote such a factor by $\FF_{r_q,s_q}$, where $r_q=j_i$ and $s_q=k_l+1$ are respectively the minimum and the maximum index plus one of the generators $\F_i$ in $\Bf_q$. For example, if $\Bf=\F_{2,4}\F_{4,5}\F_{8,8}\F_{9,9}$, then:\[\Bf_1=\F_{2,4}\F_{4,5},\quad\Bf_2=\F_{8,8}\F_{9,9}, \quad\text{ and}\quad\Bf=\FF_{2,6}\cdot\FF_{8,10}.\] Since in $\FF_{r_q,s_q}$ all $\F_i$ are present such that $r_q\le i<s_q$, to $\Bf$ there corresponds a ramified partition $(I,R)$, where\[R=\{\{r_1,r_1+1,\dots,s_1,r'_1,r'_1+1,\dots,s'_1\},\dots,\{r_m,r_m+1,\dots,s_m,r'_m,r'_m+1,\dots,s'_m \}\}\] For instance for $\Bf$ in the example above,\[R=\{\{1,1'\}\{2,3,4,5,6,2',3',4',5',6'\},\{7,7'\}\{8,9,10,8',9',10'\}\}.\] The corresponding diagram and partition diagram are shown in Figure \ref{d17}.

\begin{figure}[H]
\[\vcenter{\hbox{
\includegraphics{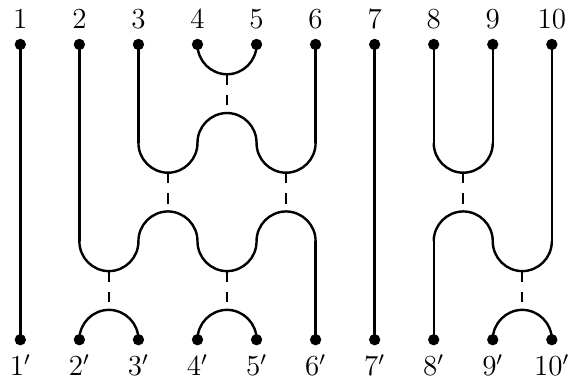}
}}
\qquad
\vcenter{\hbox{
\includegraphics{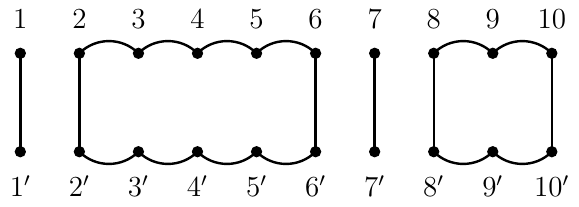}
}}\]
\caption{For $n=10$, diagrams of $\Bf$ and of the corresponding partition $R$.}
\label{d17}
\end{figure}

Observe that each block of the partition $R$ either contains a unique pair $\{i,i'\}$ or it contains all pairs $i,i'$ for $r_q\le i\le s_q$. Therefore  $R$ is  boxed.  Evidently, $N(\Bf)=\sum_q N(\Bf_q)$. For the identity, $N({\bf 1})=0$, and the number of blocks is $n$; every time a new generator $\F_j$ is introduced in $\Bf$, two blocks merge and then their number decreases by one.
\end{proof}

\begin{lemma}
Let $\F\E\in\tJ_n$  as in Proposition \ref{proEF}, involving $m$ different generators and $\psi(\F\E)=(I,Y)$, then   $Y$ is  boxed and has   $n-m$ blocks.
\end{lemma}
\begin{proof}
If $\F$ involves $k<m$ distinct generators, then $\E=\E_{i_1}\dots\E_{i_{m-k}}$, and $Y=R\ast E_{i_1}\ast\dots\ast E_{i_{m-k}}$, where $R$ is the boxed partition with $k$ blocks, given by Lemma \ref{lemmapart}. At each multiplication by $E_{i_j}$ two blocks merge, and the resulting partition is still boxed with a number of blocks diminished by one.
\end{proof}

\begin{definition}
An element $I\in J_n$ is said to be {\it separable}, if it is the union of two planar partitions, the first one  of  $[m]\cup[m']$  and  the  other of $[n-m]\cup[n'-m']$, otherwise  it is inseparable. $I$ is said $k$--separable if it is the union of $k>1$ inseparable planar partitions.
\end{definition}

\begin{example}
$\{\{1,4\},\{2,3\},\{1',2'\},\{3',4'\}\}$ is inseparable, and\\$\{\{1,3'\},\{2,3\},\{1',2'\},\{4,4'\},\{5,6\}\{5',6'\}\}$ is 3--separable into $\{\{1,3'\},\{2,3\},\{1',2'\} \}$, $\{\{4,4'\}\}$ and $\{\{5,6\}\{5',6'\}\}$.
\end{example}

\begin{proposition}\label{|bJ|}
The size of $\bJ_n$ is $\binom{2n-1}{n}$.
\end{proposition}
\begin{proof}
Let $B(n,j)$ be the number of ramified partitions $(I,R)$ such that $I$ is planar and $R$ is boxed and has $j$ blocks. We prove that\begin{equation}\label{Bnj}B(n,j)=\sum_{k=j}^{n}\binom{k-1}{j-1}T(n,n-k).\end{equation}If  $I$ is inseparable, then $R$, being boxed, has a unique block, since $I\preceq R$. For the same argument, each inseparable subpartition of $I$ is entirely contained in a block of $R$. So, $R$ has at most $k$ blocks if $I$ is  $k$--separable. The number of $k$--separable planar partitions is equal to the number of basic elements of $\J_n$ involving $n-k$ different generators, i.e. $T(n,n-k)$. This follows from Lemma \ref{lemmapart}, by recalling that  the normal forms for the elements of $\tJ_n$ with $\E=1$, correspond to those of $\J_n$. A boxed partition $R$ with $j$ blocks is thus obtained by boxed partitions with $k\ge j$ blocks by merging pairs of these blocks. There  are $\binom{k-1}{j-1}$ ways to do this, see Lemma \ref{dimlin}. Finally, taking the sum over all $j$ from $1$ to  $n$ of (\ref{Bnj}) we get\[\sum_{j=1}^n\sum_{k=j}^{n}\binom{k-1}{j-1}T(n,n-k)=\sum_{h=1}^{n-1}\left(\sum_{i=0}^{n-1-h}\binom{n-1-h}{i}\right)T(n,h)=\sum_{h=1}^{n-1}2^{n-1-h}T(n,h).\]The last sum is equal to $\binom{2n-1}{n}$ according to Lemma \ref{GCatalan}, since $T(n,n)=0$. The proof is concluded.
\end{proof}
  
Here the triangle $B(n,j)$ for $n\leq10$. Observe that the rows sum to $|\tJ_n|$, see Proposition \ref{dimtJ}.
\[\begin{tabular}{|c|cccccccccc|c|}\hline
\backslashbox{$n$}{$j$}&$1$&$2$&$3$&$4$&$5$&$6$&$7$&$8$&$9$&$10$&sum\\\hline
$1$&$1$&&&&&&&&&&$1$\\
$2$&$2$&$1$&&&&&&&&&$3$\\
$3$&$5$&$4$&$1$&&&&&&&&$10$\\
$4$&$14$&$14$&$6$&$1$&&&&&&&$35$\\
$5$&$42$&$48$&$27$&$8$&$1$&&&&&&$126$\\
$6$&$132$&$165$&$110$&$44$&$10$&$1$&&&&&$462$\\	
$7$&$429$&$572$&$429$&$208$&$65$&$12$&$1$&&&&$1716$\\
$8$&$1430$&$2002$&$1638$&$910$&$350$&$90$&$14$&$1$&&&$6435$\\
$9$&$4862$&$7072$&$6188$&$3808$&$1700$&$544$&$119$&$16$&$1$&&$24310$\\
$10$&$16796$&$25194$&$23256$&$15504$&$7752$&$2907$&$798$&$152$&$18$&$1$&$92378$\\
\hline\end{tabular}\]

\begin{proof}[Proof of Theorem \ref{isotj}]Recall that $\widetilde{E}_i=(1,E_i)$ and $\widetilde F_i=(H_i,E_i)$. Therefore if $\Bf\in\tJ_n$ and $\varphi(\Bf)=(I,R)$, then $I$ is planar since all $H_i$'s are planar partitions. Thus $\varphi(\tJ_n)\subseteq\bJ_n$. This fact together with Proposition \ref{|bJ|} and Theorem \ref{dimtJ}, says that the prooof is concluded.
\end{proof}

\end{document}